\documentclass[11pt]{article}
\usepackage[cm]{fullpage}

\usepackage{latexsym}
\usepackage{amsfonts}
\usepackage{amsmath}
\usepackage{amssymb}
\numberwithin{equation}{section}
\newcommand{\qed}{\hfill \ensuremath{\Box}}

\def\XXint#1#2#3{{\setbox0=\hbox{$#1{#2#3}{\int}$}
\vcenter{\hbox{$#2#3$}}\kern-.5\wd0}}

\newcommand{\ga}{\gamma}
\newcommand{\de}{\delta}

\newcommand{\dbar}{\overline{\partial}}

\newcommand{\ddt}[1]{\frac{\partial #1}{\partial t}}

\begin{document}
\newcounter{remark}
\newcounter{theor}
\setcounter{remark}{0} \setcounter{theor}{1}
\newtheorem{claim}{Claim}
\newtheorem{theorem}{Theorem}[section]
\newtheorem{proposition}{Proposition}[section]
\newtheorem{lemma}{Lemma}[section]
\newtheorem{definition}{Definition}[section]
\newtheorem{corollary}{Corollary}[section]
\newenvironment{proof}[1][Proof]{\begin{trivlist}
\item[\hskip \labelsep {\bfseries #1}]}{\end{trivlist}}
\newenvironment{remark}[1][Remark]{\addtocounter{remark}{1} \begin{trivlist}
\item[\hskip \labelsep {\bfseries #1
\thesection.\theremark}]}{\end{trivlist}}

\centerline{\bf THE K\"AHLER-RICCI FLOW ON SURFACES }
\centerline{\bf OF POSITIVE KODAIRA DIMENSION}
\bigskip

\begin{center}{\small
\begin{tabular}{ccc}

{\bf \normalsize Jian Song}  & & {\bf \normalsize Gang Tian} \\
Johns Hopkins University & & Princeton University \\
Department of Mathematics & & Department of Mathematics \\
3400 North Charles Street & & Fine Hall Washington Road
 \\
Baltimore MD 21218 & & Princeton NJ 08544 \\
 jsong@math.jhu.edu & & tian@math.princeton.edu
\end{tabular}}
\end{center}
\bigskip

{\footnotesize \tableofcontents}

\bigskip
\bigskip

\section{Introduction}

The existence of K\"ahler-Einstein metrics on a compact K\"ahler
manifold has been the subject of intensive study over the last few
decades, following Yau's solution to Calabi's conjecture (see
\cite{Ya2}, \cite{Au}, \cite{Ti1}, \cite{Ti2}). The Ricci flow,
introduced by Richard Hamilton in \cite{Ha1, Ha2}, has become one
of the most powerful tools in geometric analysis. The Ricci flow
preserves the K\"ahlerian property, so it provides a natural flow
in K\"ahler geometry, referred as the K\"ahler-Ricci flow. Using
the K\"ahler-Ricci flow, Cao \cite{Ca} gave an alternative proof
of the existence of K\"ahler-Einstein metrics on a compact
K\"ahler manifold with trivial or negative first Chern class. In
early 90's, Hamilton and Chow also used the Ricci flow to give
another proof of classical uniformization for Riemann surfaces
(see \cite{Ha2}, \cite{Ch}). Recently Perelman \cite{Pe} has made
a major breakthrough in studying the Ricci flow. The convergence
of the K\"ahler-Ricci flow on K\"ahler-Einstein Fano manifolds was
claimed by Perelman and it has been generalized to any K\"ahler
manifolds admitting a K\"ahler-Ricci soliton by the second named
author and Zhu \cite{TiZhu}. Previously, in \cite{ChTi}, Chen and
the second named author proved that the K\"ahler-Ricci flow
converges to a K\"ahler-Einstein metric if the initial metric is
of non-negative bisectional curvature.

However, most projective manifolds do not have a definite or
trivial first Chern class. It is a natural question to ask if
there exist any well-defined canonical metrics on these manifolds
or on varieties canonically associated to them, i.e., canonical
models. Tsuji \cite{Ts} used the K\"ahler-Ricci flow to prove the
existence of a canonical singular K\"ahler-Einstein metric on a
minimal projective manifold of general type. In this paper, we
propose a program of finding canonical metrics on canonical models
of projective varieties of positive Kodaira dimension. We also
carry out this program for minimal algebraic surfaces. To do it,
we will study the K\"ahler-Ricci flow starting from any K\"ahler
metrics and show its limiting behavior at time infinity.

Let $X$ be an $n$-dimensional compact K\"ahler manifold. A
K\"ahler metric can be given by its K\"ahler form $\omega$ on $X$.
In local coordinates $z_1, ..., z_n$, we can write $\omega$ as
$$\omega=\sqrt{-1}\sum_{i, j=1}^n
g_{i\bar{j}}dz_i\wedge dz_{\bar{j}},$$ where $\{g_{i\bar{j}}\}$ is
a positive definite hermitian matrix function.
Consider the K\"ahler-Ricci flow
\begin{equation}
\label{krflow1} \frac{\partial \omega(t,\cdot)}{\partial t} = -
 Ric(\omega(t,\cdot))- \omega(t,\cdot),
~~~~\omega(0,\cdot)=\omega_0,
\end{equation}
where $\omega(t,\cdot)$ is a family of K\"ahler metrics on $X$ and
$ Ric(\omega(t,\cdot))$ denotes the Ricci curvature of
$\omega(t,\cdot)$ and $\omega_0$ is a given K\"ahler metrics.
If the canonical class $K_X$ of $X$ is ample and $\omega_0$
represents $K_X$, Cao proved in \cite{Ca} that (\ref{krflow1}) has
a global solution $\omega(t,\cdot)$ for all $t>0$ and
$\omega(t,\cdot)$ converges to a K\"ahler-Einstein metric on $X$.
If $K_X$ is numerically effective (i.e. nef), Tsuji proved in \cite{Ts} under
additional assumption $[\omega_0]
> K_X$ that (\ref{krflow1}) has a global solution $\omega(t,\cdot)$.
This additional assumption was removed in \cite{TiZha}, moreover,
if $K_X$ is also big, $\omega(t,\cdot)$ converges to a singular
K\"ahler-Einstein metric with locally bounded K\"ahler potentials
as $t$ tends to $\infty$ (see \cite{TiZha}).
Our problem is to show how $\omega(t,\cdot)$ behaves at time
infinity.

If $X$ is a minimal K\"ahler surface of non-negative Kodaira dimension, then $K_X$ is numerically
effective. The Kodaira dimension $\nu (X)$ of $X$ is
equal to $0,1,2$. If $\nu(X)=0$, then a finite cover of $X$ is
either a K3 surface or a complex torus, then after an appropriate
scaling, $\omega(t,\cdot)$ converges to the unique Calabi-Yau
metric in the K\"ahler class $[\omega_0]$ (cf. \cite{Ca}). If
$\nu(X)=2$, i.e., $X$ is of general type, then $\omega(t,\cdot)$
converges to the unique K\"ahler-Einstein orbifold metric on its
canonical model as $t$ tends to $\infty$ (see \cite{TiZha}). If
$\nu(X)=1$, then $X$ is a minimal elliptic surface and does not
admit any K\"ahler-Einstein current which is smooth outside a
subvariety. Hence, one does not expect that $\omega(t,\cdot)$
converges to a smooth metric outside a subvariety of $X$.

In this paper, we study the limiting behavior of $\omega(t,\cdot)$
as $t$ tends to $\infty$ in the case that $X$ is a minimal
elliptic surface. In its sequel, we will extend our results here
to higher dimensional manifolds, that is, we will study the
limiting behavior of (\ref{krflow1}) when $X$ is a $n$-dimensional
variety of Kodaira dimension in $(0, n)$ and with numerically
positive $K_X$. Hence, our first question is to identify limiting
candidates. Since $X$ is a minimal elliptic surface, there is a
holomorphic map $f: X \mapsto \Sigma$ such that $K_X= \pi^*L$ for
some ample line bundle $L$ over the curve $\Sigma$. Let
$\Sigma_{reg}$ consist of all $s\in \Sigma$ such that $f^{-1}(s)$
is smooth and let $X_{reg}=f^{-1}(\Sigma_{reg})$. For any $s\in
\Sigma_{reg}$, $f^{-1}(s)$ is an elliptic curve. Thus the $L^2$-metric
on the moduli of elliptic curves induces a metric $\omega_{WP}$ on
$\Sigma_{reg}$. We call a metric $\omega$ canonical if it is
smooth on $\Sigma_{reg}$ and extends appropriately to $\Sigma$ and
satisfies
$$Ric(\omega) = -\omega + \omega_{WP},~~~~{\rm
on}~\Sigma_{reg}.$$
Such a metric exists and is unique in a
suitable sense.\footnote{Such canonical metrics can be also
defined for higher dimensional manifolds. We refer the readers to
section \ref{model} for more details.}
Here is our main result of this paper.
\begin{theorem}\label{main}
Let $f: X\rightarrow \Sigma$ be a minimal elliptic surface of
$\nu(X)=1$ with singular fibres $X_{s_1}=m_1 F_1$, ... ,
$X_{s_k}=m_k F_k$ with multiplicity $m_i\in \mathbf{N}$, $i=1,...,
k$. Then for any initial K\"ahler metric, the K\"ahler-Ricci flow
(\ref{krflow1}) has a global solution $\omega(t,\cdot)$ for all
time $t\in [0, \infty)$ satisfying:
\begin{enumerate}

\item $\omega(t,\cdot)$ converges to
$f^*\omega_{\infty}\in -2\pi c_1(X)$ as currents for a positive current
$\omega_{\infty}$ on $\Sigma$;

\item $\omega_\infty$ is smooth on
$\Sigma_{reg}$ and satisfies as currents on $\Sigma$
\begin{equation}
Ric(\omega_{\infty})=-\omega_{\infty}+\omega_{WP}+\sum_{i=1}^k\frac{m_k
-1}{m_k}[s_i],
\end{equation}
where $\omega_{WP}$ is the induced Weil-Petersson metric and
$[s_i]$ is the current of integration associated to the divisor
$s_i$ on $\Sigma$, in particular, $\omega_\infty$ is a generalized
K\"ahler-Einstein metric on $\Sigma_{reg}$;

\item for any compact subset $K\in X_{reg}$, there is a constant
$C_K$ such that
\begin{equation}
||\omega(t,\cdot)- f^*\omega_\infty||_{L^\infty(K)}+ e^t \sup_{s\in K}
||\omega(t,\cdot)|_{f^{-1}(s)}||_{L^\infty} \le C_K ,
\end{equation}
moreover, the scalar curvature of $\omega(t,\cdot)$ is uniformly
bounded on any compact set of $X_{reg}$.
\end{enumerate}

\end{theorem}

\begin{remark}
We conjecture that $\omega(t,\cdot)$ converges to
$f^*\omega_{\infty}$ in the Gromov-Hausdorff topology and in the
$C^\infty$-topology outside singular fibers.
\end{remark}

An elliptic surface $f:X\rightarrow \Sigma$ is an elliptic fibre
bundle if it does not admit any singular fibre.

\begin{corollary}\label{bundle}
Let $f: X\rightarrow \Sigma$ be an elliptic fibre bundle over a
curve $\Sigma$ of genus greater one. Then the K\"ahler-Ricci flow
 (\ref{krflow1})  has a global solution with any initial
K\"ahler metric. Furthermore, $\omega(t,\cdot)$ converges with
uniformly bounded scalar curvature to the pullback of the
K\"ahler-Einstein metric on $\Sigma$.
\end{corollary}

This theorem seems to be the first general convergence result on
collapsing of the K\"ahler-Ricci flow. Combining the results in
\cite{Ts, TiZha}, we give a metric classification for K\"aher
surfaces with an nef canonical line bundle by the K\"ahler-Ricci
flow.

\section{Preliminaries}

\subsection{Reduction of the K\"ahler-Ricci flow}

In this section, we will reduce the K\"ahler-Ricci flow
(\ref{krflow1}) to a parabolic equation on K\"ahler potentials for
any compact K\"ahler manifold $X$ with its canonical line bundle
$K_X\geq 0$.
Let $X$ be an $n$-dimensional compact K\"ahler manifold. A
K\"ahler metric can be given by its K\"ahler form $\omega$ on $X$.
In local coordinates $z_1, ..., z_n$, $\omega$ can be written in
the form of
$$\omega=\sqrt{-1}\sum_{i, j=1}^n
g_{i\bar{j}}dz_i\wedge dz_{\bar{j}},$$ where $\{g_{i\bar{j}}\}$ is
a positive definite hermitian matrix function. The curvature
tensor for $g$ is locally given by
$$ R_{i\bar{j}k\bar{l}}=-\frac{\partial^2 g_{i\bar{j}}}{\partial
z_k \partial z_{\bar{l}}}+\sum_{p, q=1}^n
g^{p\bar{q}}\frac{\partial g_{i\bar{q}}}{\partial
z_k}\frac{\partial g_{p\bar{j}}}{\partial z_{\bar{l}}}, ~~~~~i, j,
k=1, 2, ..., n.$$ And the Ricci curvature is given by
$$R_{i \bar{j}}=-\frac{\partial^2\log\det (g_{k\bar{l}})}
{\partial z_i\partial z_{\bar{j}}}, ~~~~~i, j=1, 2,..., n.$$ So
its Ricci curvature form is
$$Ric(\omega)=\sqrt{-1}\sum_{i, j=1}^n R_{i\bar{j}}dz_i\wedge
dz_{\bar{j}}=-\sqrt{-1}\partial\dbar\log \det(g_{k\bar{l}}).$$
Let $Ka(X)$ denote the K\"ahler cone of $X$, that is,
$$K_a(X)=\{ [\omega]\in H^{1,1}(X, \mathbf{R})~|~
[\omega]>0\}.$$
%
Suppose that $\omega(t,\cdot)$ is a solution of (\ref{krflow1}) on
$[0,T)$. Then its induced equation on K\"ahler classes in $Ka(X)$
is given by the following ordinary differential equation
\begin{equation}
\left\{
\begin{array}{rcl}
&&{\displaystyle \ddt{[\omega]}}  =  \displaystyle{
-2\pi c_1(X) - [\omega]}\\
&&[\omega]|_{t=0} = [ \omega_0] .
\end{array} \right.
\end{equation}
It follows
$$[\omega(t,\cdot)]=-2\pi c_1(X)+e^{-t}([\omega_0]+2\pi c_1(X)).$$
Now if we assume that $X$ has semi-positive canonical bundle
$K_X$, then there is a large integer $\mu$ such that any basis of
$H^0(X, \mu K_X)$ gives rise to an embedding $f$ into a projective
space. Recall the Kodaira dimension $\nu(X)$ of $X$ is defined to
be the dimension of the image of this embedding. This dimension is
in fact independent of choices of the basis and embedding.
Moreover, using this embedding, one can see easily that there is a
positive $(1,1)$-form $\chi$ such that $f^*\chi$ represents $-2\pi
c_1(X)$. Choose the reference K\"ahler metric $\omega_t$ by

\begin{equation}\label{ref metric}
\omega_t=\chi+e^{-t}(\omega_0-\chi).
\end{equation}
Here we abuse the notation by identifying $\chi$ and $f^*\chi$ for simplicity. Then the solution of (\ref{krflow1}) can be written as
$$\omega=\omega_t+\sqrt{-1}\partial\dbar \varphi.$$
We can always choose a smooth volume form $\Omega$ on $X$ such
that $Ric(\Omega)=\chi$. Then the evolution for the K\"ahler
potential $\varphi$ is given by the following initial value
problem:
\begin{equation} \label{**}
\left\{
\begin{array}{rcl}
&&\ddt{\varphi}=
\log\frac{e^{(n-\nu(X))t}(\omega_t+\sqrt{-1}\partial\dbar\varphi)^2}{\Omega
}-
\varphi\\
&&\varphi|_{t=0} = 0 . \end{array} \right.
\end{equation}
The following was proved in
\cite{TiZha}.
When $\omega_0>\chi$ and $K_X\ge 0$, it was proved by Tsuji in
\cite{Ts}. It was also proved in \cite{CaLa} under a stronger
assumption.
\begin{theorem}\label{longtime}
The K\"ahler-Ricci flow (\ref{**}) has a global solution for all
time $t\in [0, \infty)$ if $K_X$ is nef.
\end{theorem}
The evolution equation for the scalar curvature $R$ is given by
\begin{equation}\label{sca}
\ddt{R}=\Delta R+|Ric|^2+R.
\end{equation} Then the following proposition is an
immediate conclusion from the maximum principle for the parabolic
equation (\ref{sca}).
\begin{proposition} \label{upS} The scalar curvature along the K\"ahler-Ricci flow
 (\ref{krflow1}) is uniformly bounded from
below.
\end{proposition}

\subsection{Minimal surfaces with positive Kodaira
dimension} An elliptic fibration of a surface $X$ is a proper,
connected holomorphic map $f: X\rightarrow \Sigma$ from $X$ to a
curve $\Sigma$ such that the general fibre is a non-singular
elliptic curve. An elliptic surface is a surface admitting an
elliptic fibration. Any surface $X$ of $\nu(X)=1$ must be an
elliptic surface. Such an elliptic surfaces is sometimes called a
properly elliptic surface. Since we assume that $X$ is minimal,
all fibres are free of $(-1)$-curves. A very simple example is the
product of two curves, one elliptic and the other of genus $\geq
2$.

Let $f: X \rightarrow \Sigma$ be an elliptic surface. The
differential $df$ can be viewed as an injection of sheaves
$f^*(K_{\Sigma})\rightarrow \Omega^1_X$. Its cokernel
$\Omega_{X/\Sigma}$ is called the sheaf of relative differentials.
In general, $\Omega_{X/\Sigma}$ is far from being locally free. If
some fibre has a multiple component, then $df$ vanishes along this
component and $\Omega_{X/\Sigma}$ contains a torsion subsheaf with
one-dimensional support. Away from the singularities of $f$ we
have the following exact sequence
$$0\rightarrow f^*(K_{\Sigma})
\rightarrow\Omega^1_X\rightarrow \Omega_{X/\Sigma}\rightarrow 0$$
including an isomorphism between $\Omega_{X/\Sigma}$ and
$K_X\otimes f^*(K_{\Sigma}^{\vee})$. We also call the line bundle
$\Omega_{X/\Sigma}$ the dualizing sheaf of $f$ on $X$. The
following theorem is well-known (cf. \cite{Ba}).

\begin{theorem}\label{formula1}
Let $f: X \rightarrow \Sigma$ be a minimal elliptic surface such
that its multiple fibres are $X_{s_1}=m_1F_1$, ..,
$X_{s_k}=m_kF_k$. Then
\begin{equation}
K_X=f^{*}(K_{\Sigma}\otimes
(f_{*1}\mathcal{O}_X)^{\vee})\otimes\mathcal{O}_X(\sum
(m_i-1)F_i),
\end{equation}
or
\begin{equation}\label{eqnclass}
K_X=f^{*}(L\otimes\mathcal{O}_X(\sum(m_i-1)F_i),
\end{equation}
where $L$ is a line bundle of degree
$\chi(\mathcal{O}_X)-2\chi(\mathcal{O}_{\Sigma})$ on $\Sigma$.
\end{theorem}

Note that
$\deg(f_{*1}\mathcal{O}_X)^{\vee}=\deg(f_*\Omega_{X/\Sigma})\geq
0$ and the equality holds if and only if $f$ is locally trivial.
The following invariant
$$\delta(f)=\chi(\mathcal{O}_X)+\left(2g(\Sigma)
-2+\sum_{i=1}^{k}(1-\frac{1}{m_i})\right)$$
determines the
Kodaira dimension of $X$.

\begin{proposition}\textnormal{(cf. \cite{Ba})} \label{canonical_bundle_for
mula2}
Let $f: X\rightarrow \Sigma$ be a relatively minimal elliptic
fibration and $X$ be compact. Then $\nu(X)=1$ if and only if
$\delta(f)>0$.
\end{proposition}

Kodaira classified all possible singular fibres for $f$. A fibre
$X_s$ is stable if

\begin{enumerate}

\item $X_s$ is reduced, \item $X_s$ contains no $(-1)$-curves,

\item $X_s$ has only node singularities.
\end{enumerate}

The only stable singular fibres are of type $I_b$ for $b>0$,
therefore such singular fibres are particularly interesting.
Let $\mathcal{S}_1/\Gamma_1\cong\mathbf{C}$ be the period domain,
where $\mathcal{S}_1=\{ z\in \mathbf{C}~|~\textnormal{Im}z>0\}$ is
the upper half plane and $\Gamma_1=\textnormal{SL}(2,
\mathbf{Z})/\{\pm 1\}$ is the modular group acting by
$z\rightarrow \frac{az+b}{cz+d}$. The $j$-function gives an
isomorphism $\mathcal{S}_1/ \Gamma_1\rightarrow \mathbf{C}$ with
\begin{enumerate}
\item $j(z)=0$ if $z=e^{\frac{\pi}{3}\sqrt{-1}}$ modular
$\Gamma_1$, \item $j(z)=1$ if $z=\sqrt{-1}$ modular $\Gamma_1$.
\end{enumerate}
Thus any elliptic surface $f: X\rightarrow \Sigma$ gives a period
map $p:\Sigma_{reg} \rightarrow \mathcal{S}_1/ \Gamma_1$. Set $J:
\Sigma_{reg}\mapsto \mathbf{C}$ by $J(s)=j(p(s))$. For a stable
fibre $X_s$ of type $I_b$, the functional invariant $J$ has a pole
of order $b$ at $s$ and the
monodromy is given by
$\left(
\begin{array}{cc}
  1 & b \\
  0 & 1 \\
\end{array}
\right)$.
Now choose a semi-positive $(1,1)$-form $\chi\in-2\pi c_1(X)$ to
be the pullback of a K\"ahler form $\chi$ on the base $\Sigma$ and
$\chi=f^*\chi$ might vanish somewhere due to the presence of
singular fibres. Then as we discussed in the previous subsection,
one can reduce the K\"ahler-Ricci flow (\ref{krflow1}) to an
evolution equation on K\"ahler potentials. Hence, the
K\"ahler-Ricci flow (\ref{krflow1}) has a global solution and
provides a canonical way of deforming any given K\"ahler metric to
a canonical metric on the canonical model of minimal elliptic
surfaces of positive Kodaira dimension. As described in Theorem
\ref{main}, this canonical metric on $\Sigma$ satisfies the
curvature equation
$$
Ric(g_{\infty})=-g_{\infty}+g_{WP}+\sum_{i=1}^k\frac{m_i-1}{m_i}[s_i].$$
It corresponds to (\ref{eqnclass}), where the pullback of the
Weil-Petersson metric $g_{WP}$ by the period map $p$ is the
curvature of the dualizing sheaf $f_* \Omega_{X/\Sigma}$ and the
current $\sum_{i=1}^{k} \frac{m_i-1}{m_i}[s_i]$ corresponds to the
residues from the multiple fibres.


\section{A parabolic Schwarz lemma}\label{schwarz}

In this section we will establish a parabolic Schwarz lemma for
compact K\"ahler manifolds. This is inspired by \cite{Ya1}. This
will lead us to identify and estimate the collapsing on the
vertical direction for properly minimal elliptic surfaces and in
general certain fibre spaces. It also plays a key role in
estimating the scalar curvature along the K\"ahler-Ricci flow. Let
$f: X \rightarrow Y$ be a non-constant holomorphic mapping between
two K\"ahler manifolds. Suppose $\dim X=n$ and the K\"ahler metric
$\omega(t,\cdot)$ on $X$ is deformed by the following
K\"ahler-Ricci flow (\ref{krflow1}). Then we have the following
parabolic Schwarz lemma for metrics.

\begin{theorem}\label{schwarz1} If the holomorphic bisectional
curvature of $Y$ with respect to a fixed K\"ahler metric
$h_{\alpha\overline{\beta}}$ is bounded from above by a negative
constant $-K$ and the K\"ahler-Ricci flow (\ref{krflow1}) exists
for all $t\in [0, T)$, then
\begin{equation} f^*h\leq
\frac{C(t)}{K}\omega(t,\cdot)~,
\end{equation}
where $C(t)$ is a positive function in $t$ dependent on the
initial metric $\omega_0$ and $\lim_{t\rightarrow \infty} C(t)=1$
if $T=\infty$.
\end{theorem}

\begin{proof}
Choose normal coordinate systems for $g=\omega(t,\cdot)$ on $X$
and $h$ on $Y$ respectively. Let
$u=tr_{g}(h)=g^{i\overline{j}}f^{\alpha}_i
f^{\overline{\beta}}_{\overline{j}}h_{\alpha\overline{\beta}}$ and
we will calculate the evolution of $u$.
\begin{eqnarray*}\label{sch}
\Delta u &=& g^{k\overline{l}}\partial_k
\partial_{\overline{l}}(g^{i\overline{j}}f^{\alpha}_i
f^{\overline{\beta}}_{\overline{j}}h_{\alpha\overline{\beta}})\\
&=&g^{i\overline{l}}g^{k\overline{j}}R_{k\overline{l}}
f^{\alpha}_{i}
f^{\overline{\beta}}_{\overline{j}}h_{\alpha\overline{\beta}}+
g^{i\overline{j}}g^{k\overline{l}}f^{\alpha}_{i,k}f^{\overline{\beta}}_{\overline{j},\overline{l}}h_{\alpha\overline{\beta}}
-g^{i\overline{j}}g^{k\overline{l}}S_{\alpha\overline{\beta}
\ga\overline{\de}} f^{\alpha}_i
f^{\overline{\beta}}_{\overline{j}}f^{\ga}_k
f^{\overline{\de}}_{\overline{l}},
\end{eqnarray*}
where $S_{\alpha\overline{\beta} \ga\overline{\de}}$ is the
curvature tensor of $h_{\alpha\bar{\beta}}$ and the Laplacian
$\Delta$ acts on functions $\phi$ by
$$\Delta \phi=
g^{i\bar{j}}\partial_i\partial_{\bar{j}}\phi.$$
 \noindent By the definition of $u$ we have
\begin{eqnarray*}
\Delta u \geq g^{i\overline{l}}g^{k\overline{j}}R_{k\overline{l}}
f^{\alpha}_i
f^{\overline{\beta}}_{\overline{j}}h_{\alpha\overline{\beta}}+Ku^2.
\end{eqnarray*}
Now
\begin{eqnarray*}
\ddt{u}&=&-g^{i\overline{l}}g^{k\overline{j}}\ddt{g_{k\overline{l}}}
f^{\alpha}_i
f^{\overline{\beta}}_{\overline{j}} h_{\alpha\overline{\beta}}\\
&=&g^{i\overline{l}}g^{k\overline{j}}(R_{k\overline{l}}+g_{k\overline{l}}
)f^{\alpha}_i
f^{\overline{\beta}}_{\overline{j}}h_{\alpha\overline{\beta}}\\
&=&g^{i\overline{l}}g^{k\overline{j}}R_{k\overline{l}}f^{\alpha}_i
f^{\overline{\beta}}_{\overline{j}}h_{\alpha\overline{\beta}}+u,
\end{eqnarray*}
therefore
\begin{equation}
(\frac{d}{dt}-\Delta)u\leq u-Ku^2.
\end{equation}
Applying the maximal principle, we have
$$
\frac{d}{dt}u_{max}\leq u_{max}-Ku^2_{max}.
$$
Thus
$$u_{max}(t)\leq \frac{1}{K-Ce^{-t}}$$
and it proves the theorem. \qed
\end{proof}

By similar argument as in the proof of Theorem \ref{schwarz1} one
can also derive the following Schwarz lemma for volume forms with
weaker curvature bounds on the target manifold.

\begin{theorem}\label{schwarz2} Suppose $\dim X=n\geq \dim Y=\kappa$.
Let $\chi$ be the K\"ahler form on $Y$
 with respect to the K\"ahler metric $h_{\alpha\overline{\beta}}$.
If $Ric(h)\leq -K h$ for some $K>0$ and the K\"ahler-Ricci flow
(\ref{krflow1}) exists for all $t\in [0, T)$, then there exists a
constant $C>0$ dependent on the initial metric $\omega_0$ such
that
\begin{equation} \frac{\omega^{n-\kappa}\wedge
f^*\chi^{\kappa}}{\omega^n}\leq C.
\end{equation}
\end{theorem}

Suppose $2\pi c_1(X)=-[f^*\chi]$ for a K\"ahler form $\chi$ on
$Y$, i.e. $K_X$ is a semi-positive line bundle pulled back from an
ample line bundle from $Y$. Then we can remove the curvature
assumption on $Y$. From now on, for convenience, we will write
$f^*\chi$ as $\chi$. Since $c_1(X)\leq 0$, by Theorem
\ref{longtime}, the K\"aher-Ricci flow has long time existence.

\begin{theorem}\label{schwarz3}
Let $f:$ $X$ $\rightarrow$ $Y$ be a holomorphic fibration such
that $2\pi c_1(X)=-[f^*\chi]$ for some K\"ahler form $\chi$ on
$Y$. Then the K\"ahler-Ricci flow (\ref{krflow1}) exists for all
$t\in [0,\infty)$ and there exist constants $A, C>0$ such that for
all $(t, z)$,
\begin{equation}
f^*\chi(t,z)\leq C\omega(t,z)e^{A\varphi(t,z)}\max_{[0,t]\times
X}\{\log\frac{\Omega}{e^{(n-\kappa)s}\omega(s,\cdot)^{n}}e^{-A\varphi}\},
\end{equation}
where $\Omega$ is the volume form on $X$ such that
$Ric(\Omega)=f^*\chi$.

\end{theorem}
\begin{proof}
Let $u=g^{i\overline{j}}f^{\alpha}_i
f^{\overline{\beta}}_{\overline{j}}\chi_{\alpha\overline{\beta}}$
and choose normal coordinates for $g$ and $\chi$. We will
calculate the evolution for $\log u$. Note that $ \Delta\log
u=\frac{{\Delta}u}{u}-\frac{|{\nabla}u|_{{g}}^2}{u^2}$ and
\begin{equation}
 \Delta u
={g}^{i\overline{l}}{g}^{k\overline{j}}R_{k\overline{l}}
f^{\alpha}_{i} f^{\overline{\beta}}_{\overline{j}}+
{g}^{k\overline{l}}{g}^{i\overline{j}}f^{\alpha}_{i,k}
f^{\overline{\beta}}_{\overline{j},\overline{l}}\chi_{\alpha\overline{\beta
}}
-{g}^{i\overline{j}}{g}^{k\overline{l}}S_{\alpha\overline{\beta}
\ga\overline{\de}} f^{\alpha}_i
f^{\overline{\beta}}_{\overline{j}}f^{\ga}_k
f^{\overline{\de}}_{\overline{l}}.
\end{equation}
Applying the Cauchy-Schwartz inequality, we have
\begin{eqnarray*}
|{\nabla}u|_{{g}}^2&=&\sum_{i,j,k,\alpha,\beta}
f^{\alpha}_{i}f^{\overline{\beta}}_{\overline{j}}
f^{\alpha}_{i,k}f^{\overline{\beta}}_{\overline{j},\overline{k}}\\
&\leq&
\sum_{i,j,\alpha,\beta}|f^{\alpha}_{i}f^{\overline{\beta}}_{\overline{j}}|
(\sum_k|f^{\alpha}_{i,k}|^2)^{\frac{1}{2}}
(\sum_l|f^{\overline{\beta}}_{\overline{j},\overline{l}}|^2)^{\frac{1}{2}}\\
&=&
(\sum_{i,\alpha}|f^{\alpha}_{i}|(\sum_k|f^{\alpha}_{i,k}|^2)^{\frac{1}{2}})
^2\\
&\leq&
(\sum_{j,\beta}|f^{\beta}_{j}|^2)(\sum_{i,k,\alpha}|f^{\alpha}_{i,k}|^2).
\end{eqnarray*}
Let $C$ be a constant satisfying
$S_{\alpha\bar{\beta}\gamma\bar{\delta}}\leq C
\chi_{\alpha\bar{\beta}}\chi_{\gamma\bar{\delta}}$. Then we have
\begin{eqnarray*}
&&(\frac{d}{dt}-\Delta)\log u \\
&=&
\frac{1}{u}(-{g}^{k\overline{l}}{g}^{i\overline{j}}f^{\alpha}_{i,k}
f^{\overline{\beta}}_{\overline{j},\overline{l}}\chi_{\alpha\overline{\beta
}}
+{g}^{i\overline{j}}{g}^{k\overline{l}}S_{\alpha\overline{\beta}
\ga\overline{\de}} f^{\alpha}_i
f^{\overline{\beta}}_{\overline{j}}f^{\ga}_k
f^{\overline{\de}}_{\overline{l}}+
\frac{|{\nabla}u|_{{g}}^2}{u})+1\\
&\leq&
-\frac{1}{u}{g}^{i\overline{j}}{g}^{k\overline{l}}S_{\alpha\overline{\beta}
\ga\overline{\de}} f^{\alpha}_i
f^{\overline{\beta}}_{\overline{j}}f^{\ga}_k
f^{\overline{\de}}_{\overline{l}}+1\\
&\leq & -Cu+1.
\end{eqnarray*}
On the other hand,
\begin{eqnarray*}
(\frac{d}{dt}-\Delta)\varphi&=&-tr_{{\omega}}(\sqrt{-1}\partial\dbar
\varphi)+\ddt{\varphi}\\
&=&-tr_{{\omega}}({\omega}-
\omega_t)+\ddt{\varphi}\\
&=&tr_{{\omega}}(\omega_t)+\ddt{\varphi}-n.
\end{eqnarray*}
Notice that $\varphi$ is uniformly bounded from above from the
equation (\ref{**}) by the maximum principle. Combining the above
estimates we have
\begin{eqnarray*}
&&(\ddt{}-\Delta)(\log u-2A\varphi)\\
&\leq& -2(A-C)u- 2A\log\frac{e^{(n-\kappa)t}{\omega}^n}
{\Omega}+2A\varphi
+2nA+1 \\
&\leq&-Au+2A\log\frac{\Omega}{e^{(n-\kappa)t}\omega^n}+2A\varphi
+2nA+1 \\
&\leq& -Au+2A\log\frac{\Omega}{e^{(n-\kappa)t}\omega_t^n}+CA
\end{eqnarray*}
if we choose $A$ sufficiently large.

Suppose on each time interval $[0,t]$, the maximum of $\log
u-A\varphi$ is achieved at $(t_0, z_0)$, by the maximum principle
we have
$$u(t_0, z_0)\leq 2\log \frac{\Omega}{e^{(n-\kappa) t_0}\omega^n}(t_0, z_0)$$ and
$$u(t, z)\leq u(t_0, z_0)e ^{2A\varphi(z,t)-2A\varphi(t_0, z_0)}.$$
This completes the proof. \qed
\end{proof}


\section{Estimates}

In this section, we prove the uniform zeroth order and second
order estimate of the potential $\varphi$ along the K\"ahler-Ricci
flow. A gradient estimate is also derived and it gives a uniform
bound of the scalar curvature. We assume that $f: X\rightarrow
\Sigma$ is a properly minimal surface over a curve $\Sigma$ with
singular fibres over $\Delta=\{s_1, ..., s_{k}\}\subset \Sigma$.
Let $X_{s_i}= f^{-1}(s_i)$ be the corresponding singular fibre for
$i=1, ..., k$ and $[S]=\sum_{i=1}^{k}[X_{s_i}]=f^*(\sum_{i=1}^k[s_i])$ be the divisor
containing all the singular fibres. We can always find a hermitian metric $h$ on the line bundle induced by $[S]$ such that $Ric(h)$ is a multiple of $\chi$.

\subsection{The zeroth order and volume estimates}

We will  derive the zeroth order estimates for $\varphi$ and
$\frac{d\varphi}{dt}$.

\begin{lemma}
Let $\varphi$ be a solution of the K\"ahler-Ricci flow (\ref{**}).
There exists a positive constant $C$ depending only on the initial
data such that $\varphi\leq C$.
\end{lemma}

\begin{proof}
This is a straightforward application of the maximum principle.
Let $\varphi_{max}(t)=\max_{ X} \varphi(t, \cdot)$. Applying the
maximum principle, we have
\begin{eqnarray*}
\ddt{\varphi_{max}}&\leq&
\log\frac{e^t\omega_t^2}{\Omega}-\varphi_{max}\\
&\leq&\log\frac{2\chi\wedge(\omega_0-\chi)+e^{-t}(\omega_0-\chi)^2}{\Omega}
-\varphi_{max} \\
&\leq& C-\varphi_{max}.
\end{eqnarray*}
This gives a uniform upper bound for $\varphi$. \qed
\end{proof}

\begin{lemma}
There exists a positive constant $C$ depending only on the initial
data such that
\begin{equation}
\ddt{\varphi}\leq C.
\end{equation}
\end{lemma}

\begin{proof}
Differentiating on both sides of (\ref{**}) we obtain
\begin{equation}
\ddt{}(\ddt{\varphi})={\Delta}\ddt{\varphi}+1-
e^{-t}tr_{{\omega}}(\omega_0-\chi)-\ddt{\varphi},
\end{equation}
where $\Delta$ is the Laplacian operator of the metric $g$. It can
be rewritten as
$$
\ddt{}(e^t\ddt{\varphi})={\Delta}(e^t\ddt{\varphi})+e^t-
tr_{{\omega}}(\omega_0-\chi),
$$
and
\begin{equation}
\ddt{}(\ddt{\varphi}+\varphi)={\Delta}(\ddt{\varphi}+\varphi)+
tr_{{\omega}}(\chi)-1.
\end{equation}
So
$$
\ddt{}(e^t\ddt{\varphi}-\ddt{\varphi}-\varphi-e^t-t)=
{\Delta}(e^t\ddt{\varphi}-\ddt{\varphi}-\varphi-e^t-t)-
tr_{{\omega}}(\omega_0).
$$
Applying the maximum principle, we have
$$
e^t\ddt{\varphi}-\ddt{\varphi}-\varphi-e^t-t\leq C'
$$
for some uniform constant $C'$ only depending on the initial data.
Hence
$$
\ddt{\varphi}\leq \frac{e^{-t}}{1-e^{-t}}\varphi+C'\leq C.
$$
\qed

\end{proof}

\begin{lemma}
There exists a positive constant $C$ depending only on the initial
data such that
\begin{equation}
|\varphi|\leq C.
\end{equation}
\end{lemma}

\begin{proof}
It suffices to derive the lower bound for $\varphi$. Consider
$u(t, z)=\max_{X }\varphi(t, \cdot)-\varphi(t,z)\geq 0$.  Fix
$\delta>0$. For any $p>1$, since both $\varphi$ and
$\frac{\partial \varphi}{\partial t}$ are bounded from above,
using equation (\ref{**}), we have
\begin{equation}\label{C^0-1}
\int_{X}e^{p\delta u}({\omega}^2-\omega_t^2)\leq
\int_{X}e^{p\delta u}{\omega}^2\leq Ce^{-t}\int_{X}e^{p\delta
u}\omega_0^2.
\end{equation}
Calculate
\begin{eqnarray}\label{C^0-2}
&&\int_{X}e^{p\delta u}({\omega}^2-\omega_t^2)\nonumber\\
&=&\sqrt{-1}\int_{X}e^{p\delta u}\partial\dbar(-u)\wedge({\omega}+\omega_t)\nonumber\\
&=&\frac{2\sqrt{-1}}{ p\delta}\int_{X}\partial e^{\frac{p}{2}\delta
u}\wedge\dbar e^{\frac{p}{2}\delta u}\wedge({\omega}+\omega_t)\nonumber\\
&\geq&\frac{2\sqrt{-1}}{p\delta}\int_{X}\partial e^{\frac{p}{2}\delta
u}\wedge\dbar e^{\frac{p}{2}\delta u}\wedge\omega_t\nonumber\\
&\geq& \frac{\sqrt{-1}C}{p\delta}e^{-t}\int_{X}\partial
e^{\frac{p}{2}\delta u}\wedge\dbar e^{\frac{p}{2}\delta
u}\wedge\omega_0.
\end{eqnarray}
Combining (\ref{C^0-1}) and (\ref{C^0-2}) we obtain
$$\int_{X} |\nabla e^{\frac{p}{2}\delta u}|^2\omega_0^2\leq Cp\int_{X}
e^{p\delta u}\omega_0^2.$$ The Sobolev inequality
$||f||^2_{L^4}\leq C||f||^2_{H^1}$ implies that for all $p>1$
$$||e^{\delta u}||^{p}_{L^{2p}}\leq C\delta p||e^{\delta u}||^{p}_{L^p}.$$
Now we can apply  Moser's iteration by successively replacing $p$
by $2^k$ and let $k\rightarrow \infty$. Then the standard argument
shows that $$||e^{\delta u}||_{L^{\infty}}\leq C ||e^{\delta u}||_{L^1}.$$  Then we only need to bound the quantity
$||e^{\delta u}||_{L^1}$. Note that
$A\omega_0-\sqrt{-1}\partial\dbar
u\geq\chi+e^{-t}(\omega_0-\chi)+\sqrt{-1}\partial\dbar \varphi>0$
if we choose $A>0$ sufficiently large. The lemma is proved if we
apply the following proposition. It is proved by the second named
named author in \cite{Ti1} based on a result in \cite{Ho}.
%
%
\begin{proposition}
There exists $\delta>0$ and $C$ depending only on $(X, \omega_0)$
such that
\begin{equation}
\int_{X}e^{-\delta \phi} \omega_0^n \leq C,
\end{equation}
for all $\phi\in C^2(X)$ satisfying
$\omega_0+\sqrt{-1}\partial\dbar\phi>0$ and $\sup_{X}\phi=0$.
\end{proposition}
This completes the proof. \qed
\end{proof}


Since $e^t\omega^2=e^{\ddt{\varphi}+\varphi}\Omega$ and
$||\varphi||_{C^0}$ is uniformly bounded, from the uniform upper
bound for $\ddt{\varphi}$ we conclude that the normalized volume
form $e^t\omega^2$ is uniformly bounded above and a lower bound
for it will also give a lower bound for $\ddt{\varphi}$.
%
%

\begin{lemma}\label{volume} There exist constants $\lambda_1>0$ and $C>0$
such that for all $(t, z)\in [0, \infty)\times X$ we have the
following volume estimate
$$\frac{1}{C}|S|_h^{2\lambda_1}\leq\frac{{e^t\omega}^2}{\Omega}\leq C.$$
Here $h$ is a fixed hermitian metric equipped on the line bundle
induced by the divisor $[S]$ such that $Ric(h)>0$ is a multiple of
$\chi$.
\end{lemma}

\begin{proof}
It suffices to prove the lower bound of the volume form
$e^t\omega^2$. Notice that
$\log\frac{e^t{\omega}^2}{\Omega}=\ddt{\varphi}+\varphi$ and hence
the evolutions for $\log\frac{e^t{\omega}^2}{\Omega}$ and
$\varphi$ are prescribed by
\begin{equation}
(\frac{\partial}{\partial
t}-\Delta)\log\frac{e^t{\omega}^2}{\Omega}=tr_{{\omega}}(\chi)-1 ~~{\rm and}
\end{equation}
\begin{equation}
(\frac{\partial}{\partial
t}-\Delta)\varphi=tr_{{\omega}}(\omega_t)+
\log\frac{e^t{\omega}^2}{\Omega}-\varphi-2.
\end{equation}
Combining the above equations, at any point $(t, z)$ we have for
$\lambda>0$
\begin{eqnarray*}
&&(\frac{\partial}{\partial
t}-\Delta)(\log\frac{e^t{\omega}^2}{\Omega}+2A\varphi-\lambda\log|S|_h^2)\\
&=&2Atr_{\omega}(\omega_t)+tr_{\omega}(\chi)-\lambda
tr_{\omega}(Ric(h))
+2A\log\frac{e^t{\omega}^2}{\Omega}-2A\varphi-(4A+1)\\
&\geq&Atr_{\omega}(\omega_t)+2A\log\frac{e^t{\omega}^2}{\Omega}+
tr_{\omega}(A\omega_t-\lambda Ric(h))-C(A+1)\\
&\geq&Atr_{\omega}(\omega_t)+2A\log\frac{e^t{\omega}^2}{\Omega}-C(A+1)
\end{eqnarray*}
if we choose $A$ sufficiently large.
Suppose on each time interval $[0, T]$, the minimum of
$\log\frac{e^t{\omega}^2}{\Omega}+2A\varphi-\lambda\log|S|_h^2$ is
achieved at $(t_0, z_0)$, then by the maximum principle at $(t_0,
z_0)$ we have
\begin{equation}
tr_{{\omega}}(\omega_t)(t_0, z_0)\leq
2\log\frac{\Omega}{e^t{\omega}^2}(t_0, z_0)+C.
\end{equation}
But for some $\lambda>0$ we have at $(t_0, z_0)$
$$C+2\log\frac{\Omega}{e^t{\omega}^2}\geq tr_{{\omega}}(\omega_t)\geq
(\frac{\omega_t^2}{{\omega^2}})^{\frac{1}{2}}\geq
(\frac{\Omega}{{e^t\omega^2}})^{\frac{1}{2}}
(\frac{\chi\wedge\omega_0}{{\Omega}})^{\frac{1}{2}}\geq
C(|S|_h^{2\lambda}\frac{\Omega}{{e^t\omega^2}})^{\frac{1}{2}}.$$
For each $\delta>0$, there is the following elementary inequality
$$\log x\leq x^{\delta}+C_{\delta} ~~~\textnormal{ for all}~
x>0.$$ It follows that at $(t_0, z_0)$, we have for some small
$\delta<\frac{1}{2}$
$$(|S|_h^{2\lambda}\frac{\Omega}{{e^t\omega^2}})^{\frac{1}{2}}
\leq C((\frac{\Omega}{e^t{\omega}^2})^{\delta}+1)$$ and by
multiplying $|S|_h^{2\delta\lambda_1}$,
$$(|S|_h^{2\lambda+4\delta\lambda_1}
\frac{\Omega}{{e^t\omega^2}})^{\frac{1}{2}} \leq
C((|S|_h^{2\lambda_1}\frac{\Omega}{e^t{\omega}^2})^{\delta}+1).$$
We have $2\lambda+4\delta\lambda_1=2\lambda_1$ if $\lambda_1$ is
chosen by $\lambda_1=\frac{\lambda}{1-2\delta}.$
Therefore $|S|_h^{2\lambda_1}\frac{\Omega}{e^{t_0}{\omega}^2}(t_0,
z_0)\leq C$ and
$$\frac{e^t{\omega}^2}{|S|_h^{2\lambda_1}\Omega}e^{\varphi}(t, z)
\geq
\frac{e^{t_0}{\omega}^2}{|S|_h^{2\lambda_1}\Omega}e^{\varphi}(t_0, z_0).$$
Both $\varphi$ and
$\frac{e^{t_0}{\omega}^2}{|S|_h^{2\lambda_1}\Omega}(t_0, z_0)$ are
uniformly bounded from below, hence the lemma is proved. \qed
\end{proof}
This also shows that there is a uniform lower bound for
$\ddt{\varphi}$ with at worse $\log$ poles near the singular
fibres.

%
%

\begin{lemma}\label{time}

There exists a constant $C>0$ such that
\begin{equation}
\ddt{\varphi}\geq C(\lambda_1\log|S|_h^2-1).
\end{equation}

\end{lemma}

\begin{proof} We only have to show $\ddt{\varphi}$ is uniformly bounded from below. This is obtained by the previous lemma and
$$
\ddt{\varphi}=\log\frac{e^t{\omega}^2}{\Omega}-\varphi. $$ \qed
\end{proof}


\subsection{A partial second order estimate}

In this section, we slightly modify the proof of the parabolic
Schwarz lemma to derive a partial second order estimate. This will
imply that along the K\"ahler-Ricci flow (\ref{krflow1}) the
metric collapses along the fiber direction exponentially fast
outside the singular fibres.

\begin{lemma}\textnormal{(The partial second order estimate)}
For any $\delta>0$ there exists a constant $C>0$ depending on
$\delta$ such that
\begin{equation}
tr_{\omega}(\chi)\leq \frac{C}{|S|_h^{2\delta}}.
\end{equation}
\end{lemma}

\begin{proof}
By Lemma \ref{time}, for any $\delta>0$
$$|S|^{2\delta}\ddt{\varphi}\leq C.$$
Let $u=g^{i\overline{j}}f^{\alpha}_i
f^{\overline{\beta}}_{\overline{j}}\chi_{\alpha\overline{\beta}}$.
Following the similar calculation in Section \ref{schwarz}, we
have
\begin{eqnarray*}
&&(\ddt{}-\Delta)(\log |S|^{2\delta}u-3A\varphi)\\
&\leq& -2Au-3A\ddt{\varphi}+\delta tr_{\omega}(Ric(h))+CA\\
&\leq& -Au-3A\ddt{\varphi}+CA
\end{eqnarray*}
for sufficiently large $A$.
Suppose on each time interval $[0,T]$, the maximum of $\log
|S|^{2\delta}u-A\varphi$ is achieved at $(t_0, z_0)$, by the
maximum principle we have
$$(|S|^{2\delta}u)(t_0, z_0)\leq -3(|S|^{2\delta}\ddt{\varphi})(t_0, z_0)+C\leq C.$$
Combining with the uniform bound of $|\varphi|$, we can
conclude that  $|S|^{2\delta}u$ is uniformly bounded and the
theorem is proved. \qed
\end{proof}

\begin{corollary}\label{vertical1}
Let $X_s$ be a non-singular fibre for any $s\in \Sigma_{reg}$.
Then along the K\"ahler-Ricci flow (\ref{krflow1}), $\omega$
decays exponentially fast on $X_s$. Furthermore if $\Delta_s$ is
the Laplacian on $X_s$ with respect to $\omega_0|_{X_s}$, then
there exist constants $\lambda_2>0$ and $C>0$ such that
\begin{equation}
-e^{-t}\leq
\Delta_s\varphi\leq\frac{Ce^{-t}}{|S|^{2\lambda_2}(s)}.
\end{equation}
\end{corollary}

\begin{proof}
Applying the partial second order estimate, we have
$$0<e^{-t}+\Delta_s\varphi=\frac{\omega
|_{X_s}}{\omega_0|_{X_s}}=\frac{\omega\wedge\chi}{\omega_0\wedge\chi}=
\frac{\omega\wedge\chi}{\omega^2}\frac{\omega^2}{\omega_0\wedge\chi}\leq
tr_{\omega}(\chi)\frac{\omega^2}{\omega_0\wedge\chi}\leq
\frac{Ce^{-t}}{|S|^{2\lambda_2}(s)}$$ for some uniform constants
$C$ and $\lambda_2$. This proves the corollary. \qed
\end{proof}

The partial second-order estimate enables us to derive the
following strong partial $C^0$-estimate.

\begin{corollary}\label{vertical2}
There exists constants $\lambda_3>0$ and $C>0$ such that for all
$s\in \Sigma_{reg}$
$$|\sup_{X_s}\varphi-\inf_{X_s}\varphi|\leq \frac{Ce^{-t}}{|S|_h^{2\lambda_3}(s)}.$$
\end{corollary}

\begin{proof}
Let $\theta(s)$ be the smooth family of standard flat metrics on
the elliptic fibres over $\Sigma_{reg}$ such that
$\int_{X_s}\theta(s)=\int_{X_s}\omega_0$ for all $s\in
\Sigma_{reg}$. Let $\Delta_{\theta(s)}$ be the Laplacian of
$\theta(s)$ on each smooth fibre $X_s$.  By Green's formula, we
have
$$\varphi-\frac{1}{\int_{X_s}\theta(s)}\int_{X_s}\varphi\theta(s)
=\int_{X_s}\Delta_{\theta(s)}\varphi(y)(G_s(x,y)+A_s)\theta(s),$$
where $G_s(\cdot,\cdot)$ is Green's function with respect to
$\theta(s)$ and $A_s=\inf_{X_s\times X_s} G_s(\cdot, \cdot)$.
Since $(X_s,\theta(s))$ is a flat torus, one can easily show that
Green's function $G_s(\cdot,\cdot)$ is uniformly bounded below by
a multiple of ${\rm Diam}^2(X_s, \theta(s))$. However the diameter
${\rm diam}(X_s, \theta(s))$ might blow up near the singular
fibres and actually there exists $\lambda>0$ such that
$${\rm diam}(X_s, \theta(s))\leq \frac{C}{|S(s)|_h^{\lambda}}.$$
 Therefore
$A_s\geq- \frac{C}{|S(s)|_h^{2\lambda}}$ for some constant $C$ and
we have on each smooth fibre $X_s$,
$$|\sup_{X_s}\varphi-\inf_{X_s}\varphi
|\leq
C\sup_{X_s}|\Delta_{\theta(s)}\varphi||S(s)|_h^{-2\lambda}.$$
But for some $\mu>0$ and $C>0$ we have
$$|\Delta_{\theta(s)}\varphi|= |\Delta_s \varphi
\frac{\omega_0|_{X_s}}{\theta(s)}|= |\Delta_s \varphi|\left.
\frac{\omega_0\wedge\chi}{\theta(s)\wedge\chi}\right|_{X_s} \leq
\frac{Ce^{-t}}{|S(s)|_h^{\mu}},$$
where the last inequality follows from Corollary \ref{vertical1}
and Lemma \ref{Fbound}. This completes the proof of the corollary.
\qed
\end{proof}


\subsection{Gradient estimates}
In this section we will adapt the arguments in \cite{ChYa} to
obtain a uniform bound for $\left |\nabla \ddt{\varphi}\right |_g$ and the
scalar curvature $R$. Let
$u=\ddt{\varphi}+\varphi=\log\frac{e^t\omega^2}{\Omega}$. The
evolution equation for $u$ is given by
\begin{equation}
\ddt{u}=\Delta u+tr_{\omega}(\chi).
\end{equation}
We will obtain a gradient estimate for $u$, which will help us
bound the scalar curvature from below. Note that $u$ is uniformly
bounded from above, so we can find a constant $A>0$ such that
$A-u\geq 1$.
%
%
%
%
\begin{theorem}\label{gradient estimate}
There exist positive integers $\lambda_4$, $\lambda_5$ and a
uniform constant $C>0$ such that
\begin{enumerate}
\item[\textnormal{(i)}] $|S|_h^{2\lambda_4}|\nabla u|^2\leq
C(A-u),$
\item[\textnormal{(ii)}] $ -|S|_h^{2\lambda_5}\Delta u\leq
C(A-u),$
\end{enumerate}
where $\nabla$ is the gradient operator with respect to the metric
$g$.
\end{theorem}

\begin{proof}
Standard computation gives the following evolution equations for
$|\nabla u|^2$ and $\Delta u$:
\begin{equation}
(\frac{\partial}{\partial t}-\Delta)|\nabla u|^2=|\nabla
u|^2+(\nabla
tr_{\omega}(\chi)\cdot\overline{\nabla}u+\overline{\nabla}
tr_{\omega}(\chi)\cdot\nabla u)-|\nabla\nabla
u|^2-|\overline{\nabla}\nabla u|^2
\end{equation}
\begin{equation}(\frac{\partial}{\partial t}-\Delta)\Delta u=\Delta
u+g^{i\overline{l}}g^{k\overline{j}}R_{k\overline{l}}u_{i\overline{j}}+\Delta
tr_{\omega}(\chi).
\end{equation}
On the other hand,
$\nabla_i\overline{\nabla}_{\overline{j}}u=-R_{i\overline{j}}-\chi_{i\overline{j}}$,
hence
$$(\frac{\partial}{\partial t}-\Delta)\Delta
u=\Delta u-|\nabla\overline{\nabla}
u|^2-g^{i\overline{l}}g^{k\overline{j}}\chi_{i\overline{j}}u_{k\overline{l}
}+\Delta tr_{\omega}(\chi).$$
We shall now prove the first inequality. Let
$$H=|S|^{2\lambda_4}_h(\frac{|\nabla u|^2}{A-u}+tr_{\omega}(\chi)).$$
The evolution equation for $H$ is given by
\begin{eqnarray*}
&&(\frac{\partial}{\partial
t}-\Delta)(|S|^{2\lambda_4}_h(\frac{|\nabla
u|^2}{A-u})\\
&=&|S|_h^{2\lambda_4}\frac{|\nabla u|^2-|\nabla \nabla
u|^2-|\overline{\nabla}\nabla u|^2+(\nabla
tr_{\omega}(\chi)\cdot\overline{\nabla}u+\overline{\nabla}
tr_{\omega}(\chi)\cdot\nabla u)}{A-u}\\
&&-\Delta|S|_h^{2\lambda_4}\frac{|\nabla
u|^2}{A-u}-[\nabla|S|_h^{2\lambda_4}(\frac{\overline{\nabla}|\nabla
u|^2}{A-u}+\frac{|\nabla
u|^2\overline{\nabla}u}{(A-u)^2})+\overline{\nabla}|S|_h^{2\lambda_4}(\frac
{\nabla|\nabla
u|^2}{A-u}+\frac{|\nabla u|^2\nabla u}{(A-u)^2})]\\
&&-|S|_h^{2\lambda_4}(\frac{\nabla|\nabla
u|^2\cdot\overline{\nabla}u}{(A-u)^2}+\frac{\overline{\nabla}|\nabla
u|^2\cdot\nabla u}{(A-u)^2})+2|S|_h^{2\lambda_4}\frac{\nabla
u|^4}{(A-u)^3}.
\end{eqnarray*}
Also
\begin{eqnarray*}
&&(\frac{\partial}{\partial t}-\Delta)|S|^{2\lambda_4}_h
tr_{\omega}(\chi)\\
&=&|S|_h^{2\lambda_4}(\frac{\partial}{\partial
t}-\Delta)tr_{\omega}(\chi)-\Delta|S|_h^{2\lambda_4}tr_{\omega}(\chi)
-(\nabla|S|_h^{2\lambda_4}\cdot\overline{\nabla}
tr_{\omega}(\chi)+\overline{\nabla}|S|_h^{2\lambda_4}\cdot\nabla
tr_{\omega}(\chi))
\end{eqnarray*}
and
\begin{eqnarray*}
\nabla H&=&(|S|_h^{2\lambda_4}\frac{\nabla|\nabla
u|^2}{A-u}-|S|_h^{2\lambda_4}\frac{|\nabla u|^2\nabla
u}{(A-u)^2})+\nabla|S|_h^{2\lambda_4}\frac{|\nabla
u|^2}{A-u}+\nabla|S|_h^{2\lambda_4}tr_{\omega}(\chi)+|S|_h^{2\lambda_4}\nabla
tr_{\omega}(\chi).
\end{eqnarray*}
Since $tr_{\omega}(\chi)\leq C$ and $|S|_h^2$ can be considered as
functions pulled back from the base, one can easily show that
$$|\nabla|S|_h^{2\lambda_4}|^2\leq
C|S|_h^{4\lambda_4-2}tr_{\omega}(\chi)$$
and
$$|\Delta
|S|_h^{2\lambda_4}|\leq C|S|_h^{2\lambda_4-2}tr_{\omega}(\chi).$$
Also note that $|S|_h^{2\delta}u$ is a bounded function on $X$ for
any $\delta>0$. Calculate
\begin{eqnarray*}
&&(\frac{\partial}{\partial t}-\Delta)H\\
&\leq&|S|_h^{2\lambda_4-2}(C(\epsilon)(\frac{|\nabla
u|^2}{A-u}+\frac{|\nabla u|^2}{(A-u)^2})+\epsilon|\nabla
tr_{\omega}(\chi)|^2)+C(\epsilon)|S|_h^{2\lambda_4-1}\frac{|\nabla
u|^3}{(A-u)^2}|\\
&&+(\frac{\partial}{\partial t}-\Delta)|S|^{2\lambda_4}_h
tr_{\omega}(\chi)+\epsilon|S|_h^{2\lambda_4}|\nabla
tr_{\omega}(\chi)|^2+C(
\epsilon)\\
&&+(1-\epsilon)(\nabla
H\cdot\frac{\overline{\nabla}u}{A-u}+\overline{\nabla}
H\cdot\frac{\nabla u}{A-u})\\
&&-\frac{1}{2}|S|_h^{2\lambda_4}\frac{(|\nabla\nabla
u|^2+|\overline{\nabla}\nabla
u|^2)}{A-u}-\epsilon|S|_h^{2\lambda_4}\frac{|\nabla u|^4}{(A-u)^3}
\end{eqnarray*}
for small $\epsilon>0$. Similar calculation as in the proof of the
Schwarz lemma shows that $$(\frac{\partial}{\partial
t}-\Delta)|S|^{2\lambda_4}_h
tr_{\omega}(\chi)+\epsilon|S|_h^{2\lambda_4}|\nabla
tr_{\omega}(\chi)|^2\leq C.$$  At any point $(t_0, z_0)$ $H$
achieves its maximum  , by the maximum principle one has $\nabla
H=0$ and
$$0\leq -\frac{\epsilon}{2}|S|_h^{2\lambda_4}\frac{|\nabla
u|^4}{(A-u)^3}(t_0, z_0)+C.$$ Therefore $H(t_0, z_0)\leq C$ and
$$H\leq C.$$
Now we can prove the second inequality by making use of the first
one. Let $K=|S|_h^{2\lambda_5}\frac{A-\Delta u}{A-u}$ with $\Delta
u\leq A-1$. Then $\max K$ is uniformly bounded below from zero. By
standard calculation we have
\begin{eqnarray*}
&&(\frac{\partial}{\partial t}-\Delta)K\\
&=&|S|_h^{2\lambda_5}\frac{(|\overline{\nabla}\nabla u|^2-\Delta
u+g^{i\overline{l}}g^{k\overline{j}}\chi_{i\overline{j}}u_{k\overline{l}}-\Delta
tr_{\omega}(\chi))}{A-u}\\
&&+\Delta|S|_h^{2\lambda_5}\frac{\Delta
u}{A-u}+[(\overline{\nabla}|S|_h^{2\lambda_5})\cdot\overline{\nabla}(\frac{
\Delta
u}{A-u})+(\overline{\nabla}|S|_h^{2\lambda_5})\cdot\nabla(\frac{\Delta
u}{A-u})]\\
&&+|S|_h^{2\lambda_5}\frac{\nabla(\Delta
u)\cdot\overline{\nabla}u+\overline{\nabla}(\Delta u)\cdot\nabla
u}{(A-u)^2}+2|S|_h^{2\lambda_5}\frac{\Delta u|\nabla
u|^2}{(A-u)^3}+(\frac{\partial}{\partial
t}-\Delta)|S|_h^{2\lambda_5}\frac{A}{A-u}.
\end{eqnarray*}
and
\begin{eqnarray*}
&&(\frac{\partial}{\partial t}-\Delta)(K+3H)\\
&\leq&
(\nabla(K+3H)\cdot(\frac{\overline{\nabla}u}{A-u}+\frac{\overline{\nabla}|S
|_h^{2\lambda_5}}{|S|_h^{2\lambda_5}})
+\overline{\nabla}(K+3H)\cdot(\frac{\nabla
u}{A-u}+\frac{\nabla|S|_h^{2\lambda_5}}{|S|_h^{2\lambda_5}}))\\
&&+|S|_h^{2\lambda_5}\frac{(|\overline{\nabla}\nabla u|^2-\Delta
u+g^{i\overline{l}}g^{k\overline{j}}\chi_{i\overline{j}}u_{k\overline{l}}-\Delta
tr_{\omega}(\chi))}{A-u}-C|S|_h^{2\lambda_5-2}\frac{\Delta
u}{A-u}\\
&&-2|S|_h^{2\lambda_5}\frac{|\nabla\nabla
u|^2+|\overline{\nabla}\nabla
u|^2}{A-u}+C\\
&\leq&(\nabla(K+3H)\cdot(\frac{\overline{\nabla}u}{A-u}+\frac{\overline{\nabla}|S|_h^{2\lambda_5}}{|S|_h^{2\lambda_5}})
+\overline{\nabla}(K+3H)\cdot(\frac{\nabla
u}{A-u}+\frac{\nabla|S|_h^{2\lambda_5}}{|S|_h^{2\lambda_5}}))\\
&&-C|S|_h^{2\lambda_5-2}\frac{\Delta
u}{A-u}-\frac{3}{2}|S|_h^{2\lambda_5}\frac{|\nabla\nabla
u|^2+|\overline{\nabla}\nabla u|^2}{A-u}.
\end{eqnarray*}
Here we make use of the fact that
\begin{eqnarray*}
-\Delta tr_{\omega}(\chi)&\leq&
R_{i\overline{j}}\chi_{i\overline{j}}+Ctr^2_{\omega}(\chi)\\
&=&-u_{i\overline{j}}\chi_{i\overline{j}}-\chi_{i\overline{j}}^2+C\\
&\leq&C(|\overline{\nabla}\nabla u|+1).
\end{eqnarray*}
At any point $(t_0, z_0)$ where $H$ achieves its maximum, by the
maximum principle one has $\nabla(K+3H)=0$ and
$$|S|_h^{2\lambda_5}\frac{|\nabla\nabla
u|^2+|\overline{\nabla}\nabla u|^2}{A-u}(t_0, z_0)\leq C.$$ Hence
$K( t_0, z_0)\leq C$. Since $H$ is always bounded, one has $K(
t, z)\leq C$ for any $(t, z)$. \qed
\end{proof}


By the volume estimate we have the following immediate corollary.
\begin{corollary}\label{gradient2} For any $\delta>0$, there exists
$C>0$ depending on $\delta$ such that

\begin{enumerate}

\item $|S|_h^{2\lambda_4+\delta}|\nabla u|^2\leq C,$

\item $ -|S|_h^{2\lambda_5+\delta}\Delta u\leq C.$

\end{enumerate}
\end{corollary}


Now we are in the position to prove a uniform bound for the scalar
curvature. The following corollary tells that the K\"ahler-Ricci
will collapse with bounded scalar curvature away from the singular
fibres.

\begin{corollary}\label{gradient3}
Along the K\"aher-Ricci flow (\ref{krflow1}) the scalar curvature
$R$ is uniformly bounded on any compact set of $X_{reg}$. More
precisely, there exist constants $\lambda_6>0$ and $C$ such that
\begin{equation}
  -C \leq R \leq \frac{C}{|S|_h^{2\lambda_6}}.
\end{equation}
\end{corollary}

\begin{proof} It suffices to give an upper bound for $R$ by Proposition \ref{upS}.
Notice that
$R_{i\overline{j}}=-u_{i\overline{j}}-\chi_{i\overline{j}}$ and
then
$$R=-\Delta u-tr_{\omega}(\chi).$$
By Corollary \ref{gradient2} and the partial second order estimate, we
have
$$R\leq \frac{C}{|S|_h^{2\lambda_6}}.$$  \qed

\end{proof}

It will be interesting to know if sectional curvatures are
uniformly bounded on any compact set of $X_{reg}$. It is not
expected to be true for higher dimension. For example, we can
choose $X=X_1\times X_2$ where $X_1$ is a Calabi-Yau manifold and
$X_2$ is a compact K\"ahler manifold of $c_1(X_2)<0$. We can also
choose the initial metric $\omega_0(x_1,
x_2)=\omega_1(x_1)+\omega_2(x_2)$ where $Ric(\omega_1)=0$ and
$Ric(\omega_2)=-\omega_2$. Then along the K\"ahler-Ricci flow
(\ref{krflow1}), the solution $\omega(t,\cdot)$ is given by
$$\omega(t,x_1,x_2)=e^{-t}\omega_1(x_1)+\omega_2(x_2).$$
The bisectional curvature of $\omega_t$ will blow up along time if
the bisectional curvature of $\omega_1$ on $X_1$ does not vanish.

\subsection{The second order estimates}
In this section, we prove a second order estimate for the
potential $\varphi$ along the K\"ahler-Ricci flow. First we will
prove a formula which allows us to commute the $\partial\dbar$
operator and the push-forward operator for smooth functions on
$X$. Integrating along each fibre with respect to the initial
metric $\omega_0$, we get a function on $\Sigma$:
$$\overline{\varphi}=\frac{1}{\textnormal{vol}(X_s)}\int_{X_s}\varphi
\omega_0.$$ This can be considered as a push-forward of $\varphi$.

\begin{lemma}
Let $\varphi$ be a smooth function defined on $X$, then we have
\begin{equation}\label{commute}
\partial\dbar\int_{X_s}\varphi\omega_0=\int_{X_s}
\partial\dbar\varphi\wedge\omega_0.
\end{equation}
\end{lemma}

\begin{proof} It suffices to prove that the push forward and $\partial\dbar
$ are commutative. Let $\pi: \mathcal{M}\rightarrow B$ be an
analytic deformation of a complex manifold $M_0=\pi^{-1}(0)$.
Choose a sufficiently small neighborhood $\Delta\subset B$ such
that $M_{\Delta}=\pi^{-1}(\Delta)=\cup \Delta \times U_i$ with
local coordinates $(z^i_1,...,z^i_n, t)$, where $z^i$ is the
coordinate on $U_i$ and $t$ on $\Delta$. Now choose any test
function $\zeta$ on $B$ with ${\rm supp} \zeta\subset \Delta$ and
a partition of unity $\rho_i$ with ${\rm supp}\rho_i\subset
\Delta\times U_i$. Let $\varphi_i=\rho_i\varphi$. We calculate
\begin{eqnarray*}
\int_{M_{\Delta}} \zeta \partial\dbar\varphi\wedge\omega
&=&\sum_i\int_{\Delta\times U_i}\partial\dbar \zeta\wedge
\varphi_i\omega=\int_{\Delta}\partial\dbar
\zeta(\sum_i\int_{U_i}\varphi_i\omega)\\
&=&\int_{\Delta}\partial\dbar \zeta(\int_{M_t}\varphi\omega)
=\int_{\Delta}f\partial\dbar\int_{M_t}\varphi\omega.
\end{eqnarray*}
On the other hand side we have
\begin{eqnarray*}
\int_{M_{\Delta}} \zeta\partial\dbar\varphi\wedge\omega
&=&\int_{M_{\Delta}}\zeta\sum_i(\partial\dbar
\varphi_i\wedge\omega)
=\sum_i\int_{\Delta\times U_i}\zeta\partial\dbar\varphi_i\omega\\
&=&\int_{\Delta}\zeta(\sum_i\int_{U_i}\partial\dbar\varphi_i\wedge\omega)
=\int_{\Delta}\zeta(\int_{M_t}\partial\dbar\varphi\wedge\omega).
\end{eqnarray*}
Therefore
$$\int_{\Delta}\zeta\partial\dbar\int_{M_t}\varphi\omega
=\int_{\Delta}\zeta(\int_{M_t}\partial\dbar\varphi\wedge\omega)$$
for any testing function $f$ and hence
$$\partial\dbar\int_{M_t}\varphi\omega=
\int_{M_t}\partial\dbar\varphi\wedge\omega.$$ \qed
\end{proof}

\begin{lemma}There exists a constant $C>0$ such that
\begin{equation}(\frac{\partial}{\partial t}-\Delta)\log
tr_{\omega_0}({\omega})\leq
C(tr_{{\omega}}(\omega_0)+1).\end{equation}
\end{lemma}

\begin{proof}
Choose a normal coordinate system for $g_0$ such that ${g}$ is
diagonalized. By straightforward calculation we have
\begin{equation}
(\ddt{}-\Delta)tr_{\omega_0}({\omega})\leq
-tr_{\omega_0}(\omega)-g^{i\overline{i}}g^{k\overline{k}}g_{i\overline{j},
k}g_{j\overline{i}, \overline{k}}+C
tr_{\omega_0}(\omega)tr_{\omega}(\omega_0).
\end{equation}
It can also be shown that
\begin{eqnarray}
|{\bigtriangledown} tr_{\omega_0}({\omega})|^2
&=&\sum_{i,j,k}{g}^{k\overline{k}}
{g}_{i\overline{i},k}{g}_{j\overline{j},k}\nonumber\\
&\leq&\sum_{i,j}
(\sum_k{g}^{k}{g}^{k\overline{k}}|{g}_{j\overline{j},k}|^2)^
{\frac{1}{2}}\nonumber\\
&\leq&(\sum_i(\sum_k{g}^{k\overline{k}}|{g}_{i\overline{i},k}|^2)^{\frac{1}
{2}})^2\nonumber\\
&\leq&(\sum_i({g}_{i\overline{i}})^{\frac{1}{2}}
(\sum_k{g}^{i\overline{i}}{g}^{k\overline{k}}|{g}_{i\overline{i},k}|^2)^{\frac{1}{2}})^2\nonumber\\
&\leq&tr_{\omega_0}({\omega})\sum_{k,i}{g}^{i\overline{i}}{g}^{k\overline{k}}
|{g}_{k\overline{i},i}|^2\nonumber\\
&\leq&tr_{\omega_0}({\omega}){g}^{i\overline{i}}{g}^{k\overline{k}}
{g}_{i\overline{k},j}{g}_{k\overline{i},\overline{j}}.
\end{eqnarray}
Combined with the above inequalities, the lemma follows by
calculating $(\ddt{}-\Delta)\log tr_{\omega_0}({\omega})$. \qed
\end{proof}

\begin{lemma}
\begin{equation}
\Delta(e^t(\varphi-\overline{\varphi}))\leq -
tr_{{\omega}}(\omega_0)+\frac{1}{{\rm
Vol}(X_s)}tr_{\omega}(\int_{X_s}\omega_0^2)+2e^t.
\end{equation}
\end{lemma}

\begin{proof}
Applying equation (\ref{commute}), we have
\begin{eqnarray*}
\Delta(\varphi-\overline{\varphi})
&=&tr_{\omega}((\omega-\omega_t)
-tr_{\omega}(\frac{1}{{\rm Vol}(X_s)}\int_{X_s}\partial\dbar\varphi\wedge\omega_0)\\
&=&2-tr_{\omega}(\omega_t)-\frac{1}{{\rm Vol}(X_s)}tr_{\omega}(\int_{X_s}
\omega\wedge\omega_0-\int_{X_s}\omega_t\wedge\omega_0)\\
&\leq&-e^{-t}tr_{\omega}(\omega_0)+\frac{e^{-t}}{{\rm
Vol}(X_s)}tr_{\omega}(\int_{X_s}\omega_0^2)+2.
\end{eqnarray*}
\qed
\end{proof}


\begin{theorem}
\textnormal{(Second order estimates)} \label{secondorderestimate}
There exist constants $\lambda_7$, $A$ and $C>0$ such that
\begin{equation}tr_{\omega_0}(\omega)(t, z)\leq C
e^{Ae^t(\varphi-\overline{\varphi})(t,
z)-\frac{A}{|S|_h^{2\lambda_7} (t, z)}
\inf_{X\times[0,T]}(|S|_h^{2\lambda_7}e^s(\varphi-\overline{\varphi})))}+C.
\end{equation}
\end{theorem}

\begin{proof}
Put $H=|S|_h^{2\lambda_7}(\log
tr_{\omega_0}(\omega)-Ae^t(\varphi-\overline{\varphi}))$. We will
apply the maximum principle on the evolution of $H$. There exists
a constant $C>0$ such that
\begin{eqnarray*}
&&(\frac{\partial}{\partial t}-\Delta)H\\
&=&|S|_h^{2\lambda_7}(\frac{\partial}{\partial t}-\Delta)(\log
tr_{\omega_0}(\omega)-Ae^t(\varphi-\overline{\varphi}))-(\nabla
H\cdot\frac{\overline{\nabla}|S|_h^{2\lambda_7}}{|S|_h^{2\lambda_7}}+\overline{\nabla}
H\cdot\frac{\nabla|S|_h^{2\lambda_7}}{|S|_h^{2\lambda_7}})\\
&&+\frac{|\nabla|S|_h^{2\lambda_7}|^2}{|S|_h^{2\lambda_7}}(\log
tr_{\omega_0}(\omega)-Ae^t(\varphi-\overline{\varphi}))-\Delta
(|S|_h^{2\lambda_7})(\log
tr_{\omega_0}(\omega)-Ae^t(\varphi-\overline{\varphi}))\\
&\leq&C|S|_h^{2\lambda_7}tr_{\omega}(\omega_0)-A|S|_h^{2\lambda_7}
tr_{\omega}(\omega_0)-(\nabla
H\cdot\frac{\overline{\nabla}|S|_h^{2\lambda_7}}{|S|_h^{2\lambda_7}}+\overline{\nabla}
H\cdot\frac{\nabla|S|_h^{2\lambda_7}}{|S|_h^{2\lambda_7}})\\
&&+|S|_h^{2\lambda_7}[-Ae^t(\varphi-\overline{\varphi})
-Ae^t\ddt{(\varphi-\overline{\varphi})}+\frac{1}{{\rm Vol}(X_s)}
tr_{\omega}(\int_{X_s}\omega_0^2)]+Ce^t\\
&&+\frac{|\nabla|S|_h^{2\lambda_7}|^2}{|S|_h^{2\lambda_7}}(\log
tr_{\omega_0}(\omega)-Ae^t(\varphi-\overline{\varphi}))-\Delta
(|S|_h^{2\lambda_7})(\log
tr_{\omega_0}(\omega)-Ae^t(\varphi-\overline{\varphi})).
\end{eqnarray*}
Notice that
$$\frac{1}{{\rm
Vol}(X_s)}tr_{\omega}(\int_{X_s}\omega_0^2)\leq
tr_{\omega}(\chi)\sup_{X_s}(\frac{\omega_0^2}{\omega_{SF}\wedge
\chi}).$$
If we assume Lemma \ref{Fbound},  then
$e^t|S|_h^{2\lambda_7}(\varphi-\overline{\varphi})$,
$|S|_h^{2\lambda_7}\ddt{(\varphi-\overline{\varphi})}$ and
$|S|_h^{2\lambda_7} tr_{\omega}(\int_{X_s}\omega_0^2)$ are
uniformly bounded if $\lambda_7$ is chosen to be sufficiently
large. Also we have
$$\Delta|S|_h^{2\lambda_7}\leq
C|S|_h^{2\lambda_7-2}tr_{\omega}(\chi)$$
and
$$\frac{|\nabla|S|_h^{2\lambda_7}|^2}{|S|_h^{2\lambda_7}} \leq
C|S|_h^{2\lambda_7-2}tr_{\omega}(\chi)$$ for a uniform constant
$C>0$. Therefore we have
\begin{eqnarray*}
&&(\frac{\partial}{\partial t}-\Delta)H\\
&\leq&C|S|_h^{2\lambda_7}tr_{\omega}(\omega_0)-A|S|_h^{2\lambda_7}
tr_{\omega}(\omega_0)-(\nabla
H\cdot\frac{\overline{\nabla}|S|_h^{2\lambda_7}}{|S|_h^{2\lambda_7}}+\overline{\nabla}
H\cdot\frac{\nabla|S|_h^{2\lambda_7}}{|S|_h^{2\lambda_7}})+Ce^t.
\end{eqnarray*}
Assume $H$ achieves its maximum at $(t_0, z_0)$ on $[0, T]\times
X$. Applying the maximum principle, we have $\nabla H(t_0, z_0)=0$
and then
$$ \{|S|_h^{2\lambda_7}tr_{\omega}(\omega_0)\}(t_0, z_0)\leq C e^{t_0}.$$
This implies
$$\{|S|_h^{2\lambda_7}tr_{\omega_0}(\omega)\}(t_0, z_0)\leq C.$$ The theorem is proved then by
comparing $H$ at any point $(t, z)\in [0, T]\times X$ and $(t_0,
z_0)$. \qed

\end{proof}

\section{Generalized K\"ahler-Einstein metrics and the K\"ahler-Ricci flow}
\label{model}

\subsection{Limiting metrics on canonical models and Weil-Petersson metrics}

Let $X$ be a smooth projective manifold of $\dim X=n$. Suppose
that $\mu K_X$ is base point free for $\mu>>1$ and $\nu(X)=\kappa$
with $0<\kappa\leq n$. Then
$$|\mu K_X|: X\rightarrow X_{\mu}\subset\mathbf{CP}^{N_{\mu}}$$
is a holomorphic map for $\mu>>1$. Fix such $\mu$, we have a
holomorphic fibration $f: X\rightarrow X_{can}$ such that $\mu
K_X=f^*\mathcal{O}(1)$, where $X_{can}$ is the canonical model of
$X$ and coincides with $X_\mu$ for $\mu$ sufficiently large. The
abundance conjecture claims that there is such a holomorphic map
$f$ whenever $K_X$ is nef. If $K_X$ is also big, it was proved by
Kawamata. If $\kappa=n$, $X$ is a minimal model of general type
with $K_X$ big and nef., the K\"ahler-Ricci flow deforms any
K\"ahler metric onto a unique singular K\"ahler-Einstein metric on
$X$ (see \cite{Ts}, \cite{TiZha}). If $0<\kappa<n$, then for a
generic fibre $X_s$, one has $\mu K_{X_s}=f^*\mathcal{O}(1)$ for
some $\mu\in \mathbf{N}$, thus $\mu K_{X_s}$ is trivial and $X_s$
is a Calabi-Yau. We can choose $\chi$ to be a multiple of the
Fubini-Study metric of $\mathbf{CP}^{N_{\mu}}$ restricted on
$X_{can}$ such that $f^*\chi \in -2\pi c_1(X)$. Notice that
$f^*\chi$ is a smooth semi-positive $(1,1)$-form on $X$. For
simplicity, we sometimes denote it by $\chi$. Denote by
$X^0_{can}$ the set of all smooth points $s$ of $X_{can}$ such
that $X_s=f^{-1}(s)$ is a smooth fiber. Put
$X_{reg}=f^{-1}(X^0_{can})$. Clearly, it is a smooth manifold.

\begin{lemma}\label{sf}
 For any K\"ahler class $[\omega]$ on $X$, there is a
smooth function $\psi$ such that $\omega_{SF}:=\omega+ \sqrt{-1}
\partial \bar\partial \psi$ is a closed semi-flat ($1,1$)-form in the
following sense: the restriction of $\omega_{SF}$ to each smooth
$X_s\subset X_{reg}$ is Ricci flat.
\end{lemma}

\begin{proof}
On each smooth fiber $X_s$, let $\omega_s$ be the restriction of
$\omega$ to $X_s$ and $\partial_V$ and $\dbar_V$ be the
restriction of $\partial$ and $\dbar$ to $X_s$. Then by the Hodge
theory, there is a unique function $h_s$ on $X_s$ defined by
\begin{equation}
\left\{
\begin{array}{rcl}
&&\partial_V\dbar_V h_s=-\partial_V\dbar_V\log \omega_s^{n-\kappa}\\
&&\int_{X_s}e^{h_s}\omega_s=\int_{X_s}\omega_s.
\end{array} \right.
\end{equation}
By Yau's solution to the Calabi conjecture, there is a unique
$\psi_s$ to the following Monge-Amp\`ere equation
\begin{equation}
\left\{
\begin{array}{rcl}
&&\frac{(\omega_s+\sqrt{-1}\partial_V\dbar_V\psi_s)^{n-\kappa}}
{\omega_s^{n-\kappa}}=e^{h_s}\\
&&\int_{X_s}\psi_s\omega_s^{n-\kappa}=0.
\end{array} \right.
\end{equation}
Since $f$ is holomorphic, $\psi(z,s)=\psi_s(z)$ is well-defined as
a smooth function on $X_{reg}$.
\qed
\end{proof}

\begin{remark}
The function $\psi$ extends continuously to each smooth fiber
$X_s$ even if $s$ is a singular point of $X_{can}$.
\end{remark}

By the Hodge theory, there exists a volume form on $X$ such that
$\sqrt{-1}\partial\dbar\log \Omega={\chi}$. Define
\begin{equation}
\label{vol-comp} F=\frac{\Omega}{\left(%
\begin{array}{c}
  n \\
  \kappa \\
\end{array}%
\right)\omega_{SF}^{n-\kappa}\wedge \chi^{\kappa}}.
\end{equation}
We will show properties of $F$ and examine how $F$ behaves near
singular fibers. Define $h(z,s)=h_s(z)$, then
\begin{equation}
\label{vol-comp2} F=\frac{\Omega}{\left(%
\begin{array}{c}
  n \\
  \kappa \\
\end{array}%
\right)\omega^{n-\kappa}\wedge \chi^{\kappa}} e^{-h} \ge 0.
\end{equation}
It follows that $F$ has at most poles along $X_{can}\backslash
X_{can}^0$.

\begin{lemma} \label{Fso} $F$ is the pullback of a function
on $X_{can}$. Furthermore, there
exists $\epsilon>0$ such that
\begin{equation}
F\in L^{1+\epsilon}(X_{can}).
\end{equation}
\end{lemma}

\begin{proof} Since $\chi$ is the pullback from $X_{Can}$,
we have
$$\sqrt{-1}\partial_V\dbar_V\log\Omega=
\sqrt{-1}\partial_V\dbar_V\log\omega_{SF}^{n-\kappa}\wedge
\chi^{\kappa}=0$$
on each smooth fibre $X_s$. Thus $F$ is constant
along each smooth fibre $X_s$ and so it is the pullback of a
function from $X_{can}$.
Now we prove the second statement.
\begin{eqnarray*}
\int_{X_{can}}F^{1+\epsilon}\chi^{\kappa}
&=&\frac{1}{\int_{X_s}\omega_{SF}^{n-\kappa}}\int_{X}F^{1+\epsilon}\chi^{
\kappa}\wedge\omega_{SF}^{n-\kappa}\\
&=&\frac{1}{\left(%
\begin{array}{c}
  n \\
  \kappa \\
\end{array}%
\right)\int_{X_s}\omega_{SF}^{n-\kappa}}\int_{X}F^{\epsilon}\Omega\leq
C.
\end{eqnarray*}
The last inequality holds for sufficiently small $\epsilon>0$
because $F$ can have at worse pole singularities. \qed

\end{proof}

There is a canonical hermitian metric on the push-forward of the
dualizing sheaf
$f_*(\Omega^{n-\kappa}_{X/X_{can}})=(f_{*1}\mathcal{O}_X)^{\vee}$
over $X_{can}^0$.

\begin{definition}\label{candual2}
 Let $X$ be a projective manifold of complex dimension $n$.
 Suppose its canonical line bundle $K_X$ is semi-positive and
 $0<\kappa=\nu(X)<n$. Let $X_{can}$ be the canonical model of $X$ by
 pluricanonical maps. We define a canonical hermitian
 metric $h_{can}$ on $f_*(\Omega^{n-\kappa}_{X/X_{can}}) $ in the way that
for any smooth $(n-\kappa,
 0)$-form $\eta$ on a smooth fiber $X_s$,
\begin{equation}
|\eta|^2_{h_{can}}=\frac{\eta\wedge\bar{\eta}\wedge\chi^{\kappa}}
{\omega_{SF}^{n-\kappa}\wedge\chi^{\kappa}}
=\frac{\int_{X_s}\eta\wedge\bar{\eta}}
{\int_{X_s}\omega_{SF}^{n-\kappa}}.
\end{equation}
\end{definition}

Now let us recall some facts on the Weil-Petersson metric on the
moduli space $\mathcal{M}$ of polarized Calabi-Yau manifolds of
dimension $n-\kappa$ . Let $\mathcal{X}\rightarrow \mathcal{M}$ be
a universal family of Calabi-Yau manifolds.
Let $(U; t_1, ..., t_m)$ be a local holomorphic coordinate chart
of $\mathcal{M}$, where $m=\dim {\mathcal M}$. Then each
$\frac{\partial}{\partial t_i}$ corresponds to an element
$\iota(\frac{\partial}{\partial t_i})\in H^1(\mathcal{X}_t,
T_{\mathcal{X}_t})$ through he Kodaira-Spensor map $\iota$. The
Weil-Petersson metric is defined by the $L^2$-inner product of
harmonic forms representing classes in $H^1(\mathcal{X}_t,
T_{\mathcal{X}_t})$. In the case of Calabi-Yau manifolds, we can
express it as follows: Let $\Psi$ be a nonzero holomorphic
$(n-\kappa,0)$-form on the fibre $\mathcal{X}_t$ and
$\Psi\lrcorner\iota(\frac{\partial}{\partial t_i})$ be the
contraction of $\Psi$ and $\frac{\partial}{\partial t_i}$. Then
the Weil-Petersson metric is given by
\begin{equation}
(\frac{\partial}{\partial t_i}, \frac{\partial}{\partial
\bar{t_j}})_{\omega_{WP}}=-\frac{\int_{\mathcal{X}_t}\Psi\lrcorner\iota(\frac{\partial}{\partial
t_i})\wedge\overline{\Psi\lrcorner\iota(\frac{\partial}{\partial
t_i})}}{\int_{\mathcal{X}_t}\Psi\wedge\overline{\Psi}}.
\end{equation}
One can also represent $\omega_{WP}$ as the curvature form of the
first Hodge bundle $f_*
\Omega^{n-\kappa}_{\mathcal{X}/\mathcal{M}}$. Let $\Psi$ be a
nonzero local holomorphic section of $f_*
\Omega^{n-\kappa}_{\mathcal{X}/\mathcal{M}}$ and one can define
the hermitian metric $h_{WP}$ on $f_*
\Omega^{n-\kappa}_{\mathcal{X}/\mathcal{M}}$ by
\begin{equation}
|\Psi_t|^2_{h_{WP}}=\int_{\mathcal{X}_t}\Psi_t\wedge\overline{\Psi_t}.
\end{equation}
Then the Weil-Petersson metric is given by
\begin{equation}
\omega_{WP}=Ric(h_{WP}).
\end{equation}
\begin{lemma}
\begin{equation}
Ric(h_{can})=\omega_{WP}.
\end{equation}
\end{lemma}
\begin{proof}
Let $u=\frac{\Psi\wedge\overline{\Psi}}{\omega_{SF}^{n-\kappa}}$.
Notice $ \Psi$ restricted on each fibre $\mathcal{X}_t$ is a
holomorphic $(n-\kappa,0)$-form and  $\Psi\wedge\overline{\Psi}$
is a Calabi-Yau volume form, therefore $u$ is constant along each
fibre and can be considered as the pullback of a function on
$\mathcal{M}$. Then by definition
$$\omega_{WP}=-\sqrt{-1}\partial\dbar\log\int_{\mathcal{X}_t}
u\omega_{SF}^{n-\kappa}=-\sqrt{-1}\partial\dbar\log u,$$
where the last equality makes use of the fact that
$\int_{\mathcal{X}_t} \omega_{SF}^{n-\kappa}=constant$. At the
same time
$$Ric(h_{can})=-\sqrt{-1}\partial\dbar\log\frac{\Psi\wedge\overline{\Psi}
\wedge
\chi^{\kappa}}{\omega_{SF}^{n-\kappa}\wedge\chi^{\kappa}}
=-\sqrt{-1}\partial\dbar\log u.$$
This proves the lemma. \qed
\end{proof}

\begin{definition}\label{KEmetric}{\rm ({\bf Canonical metrics on
$X_{can}$})} \noindent We define the generalized K\"ahler-Einstein
metric $\omega$ in the class of $-2\pi f_* c_1(X)$ on $X_{can}^0$
with respect to the fibration $f:X\rightarrow X_{can}$ by
\begin{equation}
Ric(\omega)=-\omega+\omega_{WP}.
\end{equation}
In general if $X\rightarrow \Sigma$ is a Calabi-Yau fibration, we
can define a generalized K\"ahler-Einstein metric $\lambda
\omega\in 2\pi c_1(\Sigma)+2\pi c_1(f_*\Omega^{n-\kappa}_{X/\Sigma})$ by
\begin{equation}
Ric(\omega)=\lambda\omega+\omega_{WP},
\end{equation}
where $\lambda=-1$, $0$, $1$.
\end{definition}

 The following theorem is the main result of this section and its
proof is essentially due to the work of Kolodziej \cite{Ko1, Ko2}.

\begin{theorem}
Suppose that $X_{can}$ is smooth {\rm (}or has at worst orbifold
singularities{\rm )}, then there is a unique solution
$\varphi_{\infty}\in {\rm PSH}(\chi)\cap C^0(X_{can})$ of the
following equation on $X_{can}$
\begin{equation}\label{limit-eqn}
(\chi+\sqrt{-1}\partial\dbar\varphi)^{\kappa} =F \,e^{\varphi}
{\chi^{\kappa}}.
\end{equation}
Furthermore, $\omega=\chi+\sqrt{-1}\partial\dbar
{\varphi_{\infty}}$ is a positive closed current on $X_{can}$. If
$\omega$ is smooth on $X_{can}^0$, then the Ricci curvature of
$\omega$ on $X_{can}^0$ is given by
\begin{equation}\label{limit-eqn1}
Ric(\omega)=-\chi_{\omega}+\omega_{WP}.
\end{equation}
\end{theorem}

In fact, these canonical metrics belong to a class of K\"ahler
metrics which generalize Calabi's extremal metrics.
Let $Y$ be a K\"ahler manifold of complex dimension $n$ together
with a fixed closed (1,1)-form $\theta$. Fix a K\"ahler class
$[\omega]$, denote by ${\cal K}_{[\omega]}$ the space of K\"ahler
metrics within the same K\"ahler class, that is, all K\"ahler
metrics of the form $\omega_\varphi=\omega+ \sqrt{-1}\partial
\bar\partial \varphi$.
One may consider the following equation:
\begin{equation} \label{gene-C-equ} \bar\partial V_\varphi = 0,
\end{equation} where $V_\varphi$ is defined by \begin{equation}
\label{gene-C-equ2} \omega_{\varphi}(V_\varphi, \cdot) = \bar
\partial (S(\omega_\varphi) - tr_{\omega_\varphi}(\theta)).
\end{equation}
Clearly, when $\theta=0$, (\ref{gene-C-equ}) is exactly the
equation for Calabi's extremal metrics. For this reason, we call a
solution of (\ref{gene-C-equ}) a generalized extremal metric.
If $Y$ does not admit any nontrivial
holomorphic vector fields, then any generalized extremal metric
$\omega_\varphi$ satisfies
$$S(\omega_\varphi) - tr_{\omega_\varphi}(\theta) = \mu,$$
where $\mu$ is the constant given by
$$\mu=\frac{n (2\pi c_1(Y) - [\theta])\cdot
[\omega]^{n-1}}{[\omega]^n}.$$ Moreover, if $2\pi c_1(Y)-
[\theta]=\lambda [\omega]$, then any such a metric satisfies
$$Ric(\omega_\varphi) = \lambda \omega_\varphi +\theta,$$
that is, $\omega_\varphi$ is a generalized K\"ahler-Einstein
metric. This can be proved by an easy application of the Hodge
theory. More interestingly, if we take $\theta$ to be the
pull-back of $\omega_{WP}$ by $f: X_{can}^0\rightarrow {\mathcal
M}_{CY}$, then we get back those generalized K\"ahler-Einstein
metrics which arise from limits of the K\"ahler-Ricci flow.
\subsection{Minimal surfaces of general type}
We start with minimal surfaces of general type. Let $X_{can}$ be
the canonical model of a minimal surface of general type from the
contraction map $f: X\rightarrow X_{can}$. Possibly, $X_{can}$ has
rational singularities of $A-D-E$-type by contracting the
$(-2)$-curves. Since $K_{can}$ is ample and $f^*K_{can}=K_X$, we
can assume the smooth closed $(1,1)$-form $\chi=f^*{\chi}\in
-2\pi c_1(X)$ and ${\chi}$ is a K\"ahler from on $X_{can}$.  It is
shown in \cite{TiZha} that the K\"ahler-Ricci flow (\ref{krflow1})
converges to the canonical metric $g_{KE}$ on $X$, which is the
pullback of the smooth orbifold K\"ahler-Einstein metric on the
canonical model $X_{can}$, although $g_{KE}$ might vanish along
those $(-2)$-curves.
\subsection{Minimal elliptic surfaces of Kodaira dimension one}
Now consider minimal elliptic surfaces. From Lemma \ref{sf}, we
know that there exists a closed semi-flat $(1,1)$-form
$\omega_{SF}$ in $[\omega_0]$.

\begin{lemma}\label{Fbound}
Let $F$ be the function on $\Sigma$ defined by
$F=\frac{\Omega}{2\omega_{SF}\wedge\chi}$ as (\ref{vol-comp}). Let
$ B\subset \Sigma$ be a small disk with center $0$ such that all
fibres $X_s$, $s\neq 0$, are smooth. There exists a constant $C>0$
such that
\begin{enumerate}

\item If $X_0$ is of  type $mI_0$, then
\begin{equation}\label{mI0}
\frac{1}{C}|s|^{-\frac{2(m-1)}{m}}< F|_B\leq
C|s|^{-\frac{2(m-1)}{m}};
\end{equation}

\item  If $X_0$ is of type $mI_b$ or $I_b^*$, $b>0$, then
\begin{equation}\label{Ib} -\frac{1}{C}|s|^{-\frac{2(m-1)}{m}}\log|s|^2\leq
F|_B \leq -C |s|^{-\frac{2(m-1)}{m}}\log|s|^2;
\end{equation}

\item If $X_0$ is of any other type, then
\begin{equation}\label{II} \frac{1}{C}\leq F|_B\leq C. \end{equation}

\end{enumerate}

\end{lemma}
\begin{proof}Let $Y$ be the fibration of $f$ over $B$.

\begin{enumerate}

\item If $X_0$ is of type $mI_0$,  we start with a fibration
$\tilde{Y}=\mathbf{C}\times {\tilde{B}}/ L$, where
$L=\mathbf{Z}+\mathbf{Z}\cdot z(w)$ is a holomorphic family of
lattices with $z$ being a holomorphic function on $\tilde{B}$
satisfying: $z(w)=z(0)+const\cdot w^{mh}$, $w$ is the coordinate
on ${\tilde{B}}$, $h\in \mathbf{N}$. The automorphism of
$\mathbf{C}\times B$ given by $(c, w)\rightarrow (c+\frac{1}{m},
e^{\frac{2\pi\sqrt{-1}}{m}}w)$ descends to $\tilde{Y}$ and
generates a group action without fixed points. We can assume that
$Y$ is the quotient of $\tilde{Y}$ by the group action. Therefore
$\omega_{SF}$ is a smooth family of Ricci-flat metrics over $B$.
Choose a local coordinate $s$ on $B$ centered around $0$, and a
covering $\{U_{\alpha}\}$ of a neighborhood $U$ of $X_0$ in $X$ by
small polydiscs. Since the function $f^*s$ vanishes to order $m$
along $X_0$, we can in each $U_{\alpha}$ choose a holomorphic
function $w_{\alpha}$ on $U_{\alpha}$ as the $m$th root of $f^*s$,
with
$$w_{\alpha}^m=f^*s$$ and on $U_{\alpha}\cap U_{\beta}$
$$w_{\alpha}=e^{\frac{2\pi\sqrt{-1}k_{\alpha\beta}}{m}}w_{\beta}$$
for $k_{\alpha\beta}\in\{0,1,...,m-1\}$. On each $U_{\alpha}$,
$ds\wedge d\overline{s}=m^2|s|^{\frac{2(m-1)}{m}}dw_{\alpha}\wedge
d\overline{w}_{\alpha}$. Then
$|s|^{\frac{2(m-1)}{m}}F$ is smooth and bounded away from zero on
$Y$. Thus (\ref{mI0}) is proved.

\item If $X_0$ is of type $I_b$, $b>0$, we can assume
$Y=\mathbf{C}\times B/ L$, where
$$L=\mathbf{Z}+\mathbf{Z}\frac{b}{2\pi\sqrt{-1}}\log s.$$
Let $\gamma_0$ be an arc passing through $0$ in $B$ and $\gamma$
be an arc on $X$ transverse to $X_0$ with $f\cdot\gamma=\gamma_0$.
We also assume that $\gamma$ does not pass through any double
point of $X_0$. $\Omega=F\omega_{SF}\wedge \chi$ is smooth and
non-degenerate and so is $\chi$ along $\gamma$. Since
$F=\frac{\Omega}{\omega_0\wedge\chi}\frac{\omega_0\wedge\chi}{\omega_{SF}
\wedge\chi}$,
it suffices to estimate the function
$\frac{\omega_0}{\omega_{SF}}|_{X_s}$ restricted to $\gamma$ near
$X_0$. Along $\gamma$, $\omega_0|_{X_s}$ pulls back to a metric
uniformly equivalent to a fixed flat metric $\omega_{\mathbf{C}}$,
so it suffices to estimate
$\frac{\omega_{\mathbf{C}}}{\omega_{SF}}|_{X_s}$. But
$$\frac{\omega_{\mathbf{C}}}{\omega_{SF}}|_{X_s}
=\frac{\int_{X_s}\omega_{\mathbf{C}}}{\int_{X_s}\omega_{SF}}
=\frac{b\log|s|}{2\pi\int_{X_s}\omega_{SF}}$$ and ${\rm
Vol}(X_s)=\int_{X_s}\omega_{SF}$ is a constant independent of $s$.
Therefore there exists a constant $C$ such that
$$-\frac{1}{C}\log|s|^2\leq F \leq -C \log|s|^2.$$

If $X_0$ is of type $mI_b$, $b>0$, we start with a fibration
$f:\tilde{Y}\rightarrow {\tilde{B}}$, where
$\tilde{Y}=\mathbf{C}\times {\tilde{B}}/ L$ and
$L=\mathbf{Z}+\mathbf{Z}\frac{m b}{2\pi\sqrt{-1}}\log w$ and $w$
is the coordinate function of ${B}$. So
$\tilde{Y_0}=C_1+C_2+...+C_{mb}$ is of type $I_{mb}$. The
automorphism $(c, w)\rightarrow (c, e^{\frac{2\pi\sqrt{-1}}{m}}w)$
of $\mathbf{C}\times \tilde{B}$ induces a fibre-preserving
automorphism of order $m$ on $\tilde{Y}$. Such an automorphism
generates a group action on $\tilde{Y}$ without fixed points and
the quotient of $\tilde{Y}$ has a singular fibre of type $mI_b$.
Then by using the same arguments for singular fibres of type
$mI_0$, we can prove (\ref{Ib}). A fibration of type $I_b^*$
($b>0$) is obtained by taking a quotient of a fibration of type
$I_{2b}$ after resolving the $A_1$-singularities. The lattices can
be locally written as
$L=s^{\frac{1}{2}}\mathbf{Z}+\mathbf{Z}s^{\frac{1}{2}}\frac{b}{2\pi\sqrt{
-1}}\log s$. Then the above argument gives the required estimate
for $F$.

 \item If $X_0$ is not of type $mI_b$, $b\geq 0$ or $I_b^*$, $b>0$,
 it must be of type  $I^*_0$, $II$,
$III$, $IV$, $IV^*$, $III^*$ or $II^*$.  Such a singular fibre is
not a stable fiber. By the table of Kodaira (cf. \cite{Ko}), the
functional invariant $J(s)$ is bounded near $0$ and $J(0)=0$ or
$1$. One can write down the table of local lattices of periods and
the periods are bounded near the singular fibre. For example, if
$X_0$ is of type $II$, then $X_0$ is a cuspidal rational curve
with $J(s)=s^{3m+1}$, $m\in \mathbf{N}\cup \{0\}$ in the local
normal representation. On each fibre $X_s$ the above fixed flat
metric $\omega_{\mathbf{C}}$ on $\mathbf{C}$ has uniformly bounded
area, therefore
$$0<\frac{1}{C}\leq\frac{\omega_{\mathbf{C}}}{\omega_{SF}}|_{X_s}
=\frac{\int_{X_s}\omega_{\mathbf{C}}}{\int_{X_s}\omega_{SF}}\leq
C.$$ The estimate is then proved by the same argument as that in
the previous case. \qed
\end{enumerate}

\end{proof}

\begin{lemma}\label{kesol}
There is a unique solution $\varphi_{\infty}$ of the following
equation on $\Sigma$
\begin{equation}\label{limit-eqn2}
\chi+\sqrt{-1}\partial\dbar\varphi= F e^{\varphi} {\chi}
\end{equation}
satisfying $\sup_{\Sigma}|\varphi|\leq C$. Furthermore, we have
$\varphi_{\infty}\in C^{0}(\Sigma_{reg})\cap
C^{\infty}(\Sigma_{reg}).$
\end{lemma}

\begin{proof} This is a corollary of Theorem \ref{Fso}, but still we give an
elementary proof for the sake of completeness. Rewrite equation
(\ref{limit-eqn2}) as
\begin{equation}\label{lap}
\Delta\varphi=Fe^{\varphi}-1,
\end{equation}
where $\Delta$ is the Laplacian operator with respect to $\chi$.
We will apply the method of continuity to find the solutions of
the following equation parameterized by $t\in[0,1]$:
\begin{equation}
\label{conti} \Delta \varphi=e^{\varphi}(\frac{2\omega_{SF}\wedge
\chi}{\Omega}+t)^{-1}-1.
\end{equation}
Obviously equation (\ref{conti}) is solvable for all $t\in(0,1]$.
 To solve for $t=0$ we need to derive
the uniform $C^0$-estimate for $\varphi_t$.  By the maximum
principle, $$\sup_{\Sigma\times (0,1]}\varphi_t
\leq\sup_{\Sigma}\log\frac{2\omega_{SF}\wedge \chi}{\Omega}\leq
C.$$ By Lemma \ref{Fbound}, $||F||_{L^p}$ is bounded for some
$p>1$, then the standard $L^p$ estimate gives 
$$||\varphi_t||_{L^p_2}\leq C(||F||_{L^p}+1)\leq C.$$ The Sobolev
embedding theorem implies $$||\varphi_t||_{L^{\infty}}\leq C $$
for $t\in (0, 1]$. With the $C^0$ estimate, we can derive the
uniform $C^k$-estimate for $\varphi_t$ by the local estimates of
the standard theory of linear elliptic PDE due to the fact that
$\Delta$ has uniformly bounded coefficients. Therefore there
exists $\varphi_{\infty}\in C^{\infty}(\Sigma_{reg})$ satisfying
equation (\ref{limit-eqn2}).

 Now we prove the uniqueness. Suppose there is another
solution $\psi$ solving (\ref{limit-eqn}) with
$\sup_{\Sigma}|\psi|\leq C$. We define
$$\varphi_{\epsilon}=\varphi+\epsilon\log|S|_h^2$$ such that $|S|_h^2\leq
1$. Then
$$
\sqrt{-1}\partial\dbar(\varphi_{\epsilon}-\psi)
=(\frac{1}{|S|^{2\epsilon}_h}e^{\varphi_{\epsilon}-\psi}-1)({\chi+\sqrt{-1}\partial\dbar\psi})-\epsilon\chi.
$$
Note that $\{\varphi_{\infty}-\psi\geq 0\}\cap (\cup_i
X_{p_i})=\phi$, then
\begin{eqnarray*}
&&-\int_{\{\varphi_{\epsilon}-\psi\geq
0\}}|\nabla_{\psi}(\varphi_{\epsilon}-\psi)|^2(\chi+\sqrt{-1}\partial\dbar\psi)\\
&=&\int_{\{\varphi_{\epsilon}-\psi\geq
0\}}(\varphi_{\epsilon}-\psi)\sqrt{-1}\partial\dbar(\varphi_{\epsilon}-\psi)\\
&=&\int_{\{\varphi_{\epsilon}-\psi\geq
0\}}(\varphi_{\epsilon}-\psi)(\frac{1}{|S|^{2\epsilon}_h}e^{\varphi_{\epsilon}-\psi}-1)({\chi+\sqrt{-1}\partial\dbar\psi})-\int_{\{\varphi_{\epsilon}-\psi\geq
0\}}\epsilon\chi\\
&\geq& -C\epsilon.
\end{eqnarray*}
Let $\epsilon\rightarrow 0$, by Fatou's lemma,
$$
\int_{\{\varphi-\psi\geq
0\}}|\nabla_{\psi}(\varphi_{\epsilon}-\psi)|^2(\chi+\sqrt{-1}\partial\dbar\psi)=0,
$$
therefore $\varphi-\psi=0$ on each path connected component of
$\{\varphi-\psi\geq 0\}.$ On the other hand, we can apply the same
argument for $\psi_{\epsilon}$ and it infers that
$$
\varphi-\psi= 0.
$$

\qed
\end{proof}

\begin{corollary} Let $f: X\rightarrow \Sigma$ be a minimal elliptic surfac
e of
$\nu(X)=1$ with singular fibres $X_{s_1}=m_1 F_1$, ... ,
$X_{s_k}=m_k F_k$ with multiplicity $m_i\in \mathbf{N}$, $i=1,...,
k$. If $\varphi_{\infty}$ is the solution in Lemma \ref{kesol},
then
 $\omega_{\infty}=\chi+\sqrt{-1}\partial\dbar \varphi_{\infty}$ is a
K\"ahler form on $\Sigma$ and Ricci curvature of $\omega_{\infty}$
is given by the following formula
\begin{equation}
Ric(\omega_{\infty})=-\omega_{\infty}+\omega_{SF}+\sum_{i=1}^k\frac{m_k
-1}{m_k}[s_i],
\end{equation}
where $\omega_{WP}$ is the induced Weil-Petersson metric and
$[s_i]$ is the current of integration associated to the divisor
$s_i$ on $\Sigma$. In particular, if $f: X\rightarrow \Sigma$ has
only singular fibres of type $mI_0$, then $\chi_{\infty}$ is a
hyperbolic cone metric on $X$ given by
\begin{equation}
Ric(\chi_{\infty})=-\omega_{\infty}.
\end{equation}

\end{corollary}
 This tells us that $\omega_{\infty}$ is canonical in the sense
that
 it is more or less a hyperbolic metric with a correction term
 $\omega_{SF}+\sum_{i=1}^k\frac{m_k-1}{m_k}[s_i]$
 inherited  from the elliptic fibration structure of $X$. Also we
 notice that the residues only come from multiple fibres.

\section{Convergence}

\subsection{Sequential convergence}

In this section we will prove a sequential convergence of the
K\"ahler-Ricci flow to a canonical metric on the base $\Sigma$.
\begin{lemma}\label{ddtbd} There exist constants $\lambda_8$ and $C$ such that
\begin{equation}
\left|\nabla_{g_0}\ddt{\varphi}\right|^2\leq
Ce^{\frac{C}{|S|_h^{2\lambda_8}}} .
\end{equation}
\end{lemma}

\begin{proof} Put $u=\ddt{\varphi}+\varphi$ and calculate
\begin{eqnarray*}
&&\frac{1}{2}\left|\nabla_{g_0}\ddt{\varphi}\right|^2
-\left|\nabla_{g_0}\varphi\right|^2\\
&\leq&
\left|\nabla_{g_0}(\ddt{\varphi}+\varphi)\right|^2=\frac{\sqrt{-1}\partial
u\wedge\dbar u\wedge\omega_0}{\omega_0^2}=\left|\nabla
u\right|^2\frac{\omega\wedge\omega_0}{\omega_0^2}\\
&\leq& \left|\nabla u\right|^2tr_{\omega_0}(\omega).
\end{eqnarray*}
The lemma follows from Theorem \ref{secondorderestimate}(the
second order estimate), Corollary \ref{vertical2} and Theorem
\ref{gradient estimate} (the gradient estimate) since
$\left|\nabla_{g_0}\varphi\right|$ can be bounded by
$tr_{\omega_0}(\omega)$. \qed
\end{proof}
We first prove a weak convergence for $\ddt{\varphi}$ by picking a
sequence.

\begin{lemma} There exists a sequence ${t_j}$ such that
$\ddt{\varphi}(t_j,\cdot)$ converges to $0$ weakly. Thus by Lemma
\ref{ddtbd}, $\ddt{\varphi}(t_j, \cdot)$  converges to $0$ in
$C^{0,\alpha}$ for any $0<\alpha<1$ on any compact set of
$X_{reg}$.
\end{lemma}

\begin{proof}
For any $T_2>T_1>0$ and $z\in X$ we have
$$\int_{T_1}^{T_2}\ddt{\varphi}(s,z)ds=\varphi( T_2, z)-\varphi( T_1, z).
$$
Since $|\varphi|_{C^0}$ is uniformly bounded on $X$, we have
$$\left|\int_{T_1}^{T_2}\ddt{\varphi}(s,z)ds\right|\leq C.$$
Choose a countable basis $\{U_{\alpha}\}_{\alpha\in \cal{A}}$ for
the topology  of $X$. Then on each $U_{\alpha}$,
$\int_{T_1}^{T_2}(\int_{U_{\alpha}}\ddt{\varphi}\omega_0^2)ds$ is
uniformly bounded independent of the choice of $T_1$ and $T_2$.
Therefore by applying the mean value theorem with $T_1, T_2,
T_2-T_1\rightarrow \infty$, we can show that by passing to a
sequence $\ddt\varphi(t_{k, \alpha}, \cdot)$
$$\lim_{k\rightarrow\infty}\int_{U_{\alpha}}\ddt{\varphi}
(t_{k,\alpha},\cdot)=0.$$ By taking the diagonal sequence $t_j$ of
$t_{k, \alpha}$, one has for each $\alpha\in \cal{A}$
$$\lim_{j\rightarrow\infty}\int_{U_{\alpha}}\ddt{\varphi}(t_j,\cdot)\omega_
0^2=0.$$
Therefore $\ddt{\varphi}(t_j,\cdot)$ converges weakly to $0$ on
$X$.\qed
\end{proof}

\begin{theorem}\label{seq-con}
There exists a sequence $\{t_j\}_{j=1}^{\infty}$ with
$t_j\rightarrow \infty$ such that $\varphi(t_j,\cdot)$ converges
to the pullback of a function $\varphi_{\infty}$ given by equation
 (\ref{limit-eqn2}) uniformly on any compact set of
$X_{reg}$ in the sense of $C^{1,1}$.
\end{theorem}

\begin{proof}
Consider any test function $\zeta\in C^{\infty}_0(\Sigma_{reg})$.
We calculate
\begin{eqnarray*}
\int_{X}\zeta e^t{\omega}^2 &=&\int_{X} \zeta\left[
2(\chi+\sqrt{-1}\partial\dbar\bar{\varphi})\wedge\omega_0+
2\sqrt{-1}e^t\partial\dbar(\varphi-\bar{\varphi})\wedge(\chi+\sqrt{-1}\partial\dbar\bar{\varphi})\right.\\
&&~~~~~~\left.+e^{-t}(\omega_0-\chi+\sqrt{-1}e^t\partial\dbar(\varphi-\bar{\varphi}))^2\right].
\end{eqnarray*}
Notice that $\int_{X}\zeta\partial
\dbar(\varphi-\bar{\varphi})\wedge(\chi+\sqrt{-1}\partial\dbar\bar{\varphi})
=\int_{X}(\varphi-\bar{\varphi})\partial
\dbar\zeta\wedge(\chi+\sqrt{-1}\partial\dbar\bar{\varphi})=0$ and
\begin{eqnarray*}
&&\int_{X}\zeta e^{-t}(\omega_0-\chi+\sqrt{-1}e^t\partial\dbar(\varphi-\bar{\varphi}
))^2\\
&=&\int_{X}\zeta\left[e^{-t}(\omega_0-\chi)^2
+2\sqrt{-1}\partial\dbar(\varphi-\bar{\varphi})\wedge(\omega_0-\chi)
+\sqrt{-1}e^{t}\partial\dbar(\varphi-\bar{\varphi})\wedge\sqrt{-1}\partial\dbar(\varphi-\bar{
\varphi})\right]\\
&=&\int_{X}\left[2(\varphi-\bar{\varphi})\sqrt{-1}\partial\dbar
\zeta\wedge(\omega_0-\chi)
+e^{t}(\varphi-\bar{\varphi})\partial\dbar
\zeta\wedge\sqrt{-1}\partial\dbar(\varphi
-\bar{\varphi})\right]+O(e^{-t})\\
&=&O(e^{-t}).
\end{eqnarray*}
Therefore
$$\int_{X}2\zeta(\chi+\sqrt{-1}\partial\dbar\bar{\varphi})\wedge\omega_0=
\int_{X}\zeta e^t{\omega}^2+O(e^{-t})= \int_{X}\zeta
e^{\varphi+\ddt{\varphi}}\Omega+O(e^{-t}).$$
Since $\chi+\sqrt{-1}\partial\dbar\bar{\varphi}$ sits on the base
$\Sigma$ and the volume along each smooth fibre given by
$\omega_{SF}$ is the same as that of $\omega_0$, we have
$$\int_{X}2\zeta(\chi+\sqrt{-1}\partial\dbar\bar{\varphi})\wedge\omega_{SF}=
\int_{X}\zeta
e^{\varphi+\ddt{\varphi}}\frac{\Omega}{\chi\wedge\omega_{SF}}\chi\wedge\omega_{SF}+O(e^{-t}).
$$
So for all $\zeta\in C^{\infty}_0(\Sigma_{reg})$ , we have
$$\lim_{t_j\rightarrow\infty}\int_{X}2\zeta(\chi+\sqrt{-1}\partial\dbar\bar{\varphi})
\wedge\omega_{SF} =\int_{X}\zeta
e^{\varphi_{\infty}}\frac{\Omega}{\chi\wedge\omega_{SF}}\chi\wedge\omega_{SF}.$$
since $\ddt{\varphi}(t_j,\cdot)\rightarrow 0$. By taking the
convergent subsequent we have
$$\frac{\chi+\sqrt{-1}\partial\dbar{\varphi_{\infty}}}{\chi}
=e^{\varphi_{\infty}}\frac{\Omega}{2\chi\wedge\omega_{SF}}.$$
Notice $\varphi_{\infty}$ is uniformly $C^0$ bounded since
$\varphi$ is uniformly bounded in $C^0(X)$ along the
K\"ahler-Ricci flow. Also by the uniqueness of the solution for
the equation above, we have proved the theorem. \qed
\end{proof}


\subsection{Uniform convergence}

In this section we will prove a uniform convergence of the
K\"ahler-Ricci flow.
Since $\varphi$ and $\varphi_{\infty}$ are both uniformly bounded
on $X$. Therefore for any $\epsilon>$, there exists
$r_{\epsilon}>0$ with $\lim_{\epsilon\rightarrow
0}r_{\epsilon}=0$, such that for any $z\in \cup_{i=1}^{\mu}
B_{r_{\epsilon}}(p_i)$ and $t>0$ we have
\begin{eqnarray*}
&&(\varphi-\varphi_{\infty}+\epsilon\log|S|_h^2)(t, z)<-1 ~~{\rm and}\\
&&(\varphi-\varphi_{\infty}-\epsilon\log|S|_h^2)(t, z)>1,
\end{eqnarray*}
where $B_{r_{\epsilon}}(p_i)$ is a geodesic tube centered at the
singular fibre $ X_{p_i}$ with radius $r_{\epsilon}$ with respect
to the metric $\omega_0$.
 Suppose
the semi-flat closed form is given by
$\omega_{SF}=\omega_0+\sqrt{-1}\partial\dbar\rho_{SF}$ and
$\rho_{SF}$ blows up near the singular fibres. We can always find
an approximation $\rho_{\epsilon}$ for $\rho_{SF}$ such that
$\rho_{\epsilon}$ is smooth on $X$ and on $X\setminus
\cup_{i=0}^{\mu}B_{r_{\epsilon}}(p_i)$
$$\rho_{\epsilon}=\rho_{SF}.$$
We also define
$\omega_{SF,\epsilon}=\omega_0+\sqrt{-1}\partial\dbar\rho_{\epsilon}$.
Now we define the twisted difference  of $\varphi$ and
$\varphi_{\infty}$ by
\begin{eqnarray}
&&\psi^{-}_{\epsilon}=\varphi-\varphi_{\infty}-e^{-t}\rho_{\epsilon}
+\epsilon\log|S|^2_{h}~~{\rm and}\\
&&\psi^{+}_{\epsilon}=\varphi-\varphi_{\infty}-e^{-t}\rho_{\epsilon}
-\epsilon\log|S|^2_{h}.
\end{eqnarray}

\begin{lemma}\label{conv}
There exists $\epsilon_0>$ such that for any
$0<\epsilon<\epsilon_0$, there exists $T_{\epsilon}>0$ such that
for and $z\in X$ and $t>T_{\epsilon}$ we have
\begin{eqnarray}
&&\psi^{-}_{\epsilon}(t, z)\leq 4\epsilon~~{\rm and}\\
&&\psi^{+}_{\epsilon}(t, z)\geq -4\epsilon.
\end{eqnarray}
\end{lemma}

\begin{proof}
 The evolution for $\psi_{\epsilon}^{-}$ is given by
\begin{equation}
\ddt{\psi^{-}_{\epsilon}}=\log
\frac{e^t(\chi_{\infty}+\epsilon\chi+e^{-t}\omega_{SF,
\epsilon}+\sqrt{-1}\partial\dbar\psi^{-}_{\epsilon})^2}{\chi_{\infty}\wedge
\omega_{SF}}-\psi^{-}_{\epsilon}+\epsilon\log|S|^2_h.
\end{equation}
Since $\rho_\epsilon$ is bounded on $X$, we can always choose
$T_1>0$ sufficiently large such that for $t>T_1$   we have
$\psi^{-}_{\epsilon}(t, z)<-\frac{1}{2}$ on
$\cup_{i=1}^{\mu}B_{r_{\epsilon}}(p_i)$ and
$e^{-t}\frac{\omega_{SF}^2}{\chi_{\infty}\wedge\omega_{SF}}\leq
\epsilon$ on $X\setminus \cup_{i=1}^{\mu}B_{r_{\epsilon}}(p_i)$.
We will discuss in two cases for $t>T_1$.

\begin{enumerate}

\item If $\psi^{-}_{\epsilon, max}(t)=\max_{
X}\psi^{-}_{\epsilon}(t, \cdot)=\psi^{-}_{\epsilon}(t,
z_{max,t})>0$ for all $t>T_1$. Then $z_{max,t}\in X\setminus
\cup_{i=1}^{\mu}B_{r_{\epsilon}}(p_i)$ for all $t>T_1$. Applying
the maximum principle at $z_{max,t}$ we have
\begin{eqnarray*}
\ddt{\psi^{-}_{\epsilon}}(t, z_{max,t})&\leq&
\{\log\frac{e^t(\chi_{\infty}+\epsilon\chi+e^{-t}\omega_{SF,\epsilon})^2}
{\chi_{\infty}\wedge\omega_{SF}}-\psi^{-}_{\epsilon}+\epsilon\log|S|^2_h\}
(t,z_{max,t})\\
&=&\{\log\frac{(\chi_{\infty}+\epsilon\chi)\wedge\omega_{SF,\epsilon}+
e^{-t}\omega_{SF,\epsilon}^2}{\chi_{\infty}\wedge\omega_{SF}}-
\psi^{-}_{\epsilon}+\epsilon\log|S|^2_h\}(t, z_{max,t})\\
&=&\{\log\frac{(\chi_{\infty}+\epsilon\chi)\wedge\omega_{SF}+
e^{-t}\omega_{SF}^2}{\chi_{\infty}\wedge\omega_{SF}}-
\psi^{-}_{\epsilon}+\epsilon\log|S|^2_h\}(t, z_{max,t})\\
&\leq&- \psi^{-}_{\epsilon}(t,
z_{max,t})+\log(1+2\epsilon)+\epsilon.
\end{eqnarray*}
Applying the maximum principle again, we have
\begin{equation}
\psi^{-}_{\epsilon}\leq 4\epsilon+O(e^{-t}).
\end{equation}

\item If there exists $t_0\geq T_1$ such that $\max_{z\in
X}\psi^{-}_{\epsilon}(t_0, z)=\psi^{-}_{\epsilon}(t_0, z_0)<0$ for
some $z_0\in X$. Assume $t_1$ is the first time when $\max_{z\in
X, t\leq t_1}\psi^{-}_{\epsilon}(t, z)=\psi^{-}_{\epsilon}(t_1,
z_1)>4\epsilon$. Then $z_1\in X\setminus
\cup_{i=1}^{\mu}B_{r_{\epsilon}}(p_i)$ and applying the maximum
principle we have
\begin{eqnarray*}
\psi^{-}_{\epsilon}(t_1, z_1)&\leq&\{
\log\frac{(\chi_{\infty}+\epsilon\chi)\wedge\omega_{SF}+
e^{-t_1}\omega_{SF}^2}{\chi_{\infty}\wedge\omega_{SF}}+
\epsilon\log|S|^2_h\}(t_1, z_1)\\
&\leq&\log(1+2\epsilon)+\epsilon< 4\epsilon.
\end{eqnarray*}
which contradicts  the assumption that $\psi^{-}_{\epsilon}(t_1,
z_1)\geq 4\epsilon$.
Hence we have
\begin{equation}
\psi^{-}_{\epsilon}\leq 4\epsilon.
\end{equation}
By the same argument we have
\begin{equation}
\psi^{+}_{\epsilon}\geq -4\epsilon.
\end{equation}
This completes the proof.
\qed
\end{enumerate}

\end{proof}

\begin{lemma}
We have the point-wise convergence of $\varphi$ on $X_{reg}$.
Namely, for any $z\in X_{reg}$ we have
\begin{equation}
\lim_{t\rightarrow \infty}\varphi(t, z)=\varphi_{\infty}(z).
\end{equation}
\end{lemma}

\begin{proof}
By lemma \ref{conv}, we have for $t>T_{\epsilon}$
\begin{equation}
\varphi_{\infty}(t, z)+\epsilon\log|S|^2_h(t, z)-4\epsilon \leq
\varphi(t, z)\leq
                \varphi_{\infty}(t, z)-\epsilon\log|S|^2_h(t, z)+4\epsilon.
\end{equation}
Then the lemma is proved by letting $\epsilon\rightarrow 0$.
\end{proof}

Since we have the uniform zeroth and second order estimates for
$\varphi$ away from the singular fibres, we derive our main
theorem.
\begin{theorem}\label{uni-con}
 $\varphi$ converges to the
pullback of a function $\varphi_{\infty}$ given by equation
(\ref{limit-eqn}) on $\Sigma$ uniformly on any compact set of
$X_{reg}$ in the sense of $C^{1,1}$.
\end{theorem}


\section{An alternative deformation and large complex structure
 limits}
Mirror symmetry and the SYZ conjecture make predictions for
Calabi-Yau manifolds with "large complex structure limit point"
(cf. \cite{SYZ}). It is believed that in the large complex
structure limit, the Ricci-flat metrics should converge in the
Gromov-Hausdorff sense to a half-dimensional sphere by collapsing
a special Lagrangian torus fibration over this sphere. This holds
trivially for elliptic curves and is proved by Gross and Wilson
(cf. \cite{GrWi}) in the case of $K3$ surfaces. The method of the
proof is to find a good approximation for the Ricci-flat metrics
near the large complex structure limit. The approximation metric
is obtained by gluing together the Oogrui-Vafa metrics near the
singular fibres and a semi-flat metric on the regular part of the
fibration. Such a limit metric of $K3$ surfaces is McLean's
metric.

In this section, we will apply a deformation for a family of
Calabi-Yau metrics and derive  Mclean's metric \cite{Mc} without
writing down an accurate approximation metric. Such a deformation
can be also done in higher dimensions. It will be interesting to
have a flow which achieves this limit. The large complex structure
limit of a $K3$ surface $\hat{X}$ can be identified as the mirror
to the large K\"ahler limit of $X$ as shown in \cite{GrWi}, so we
can fix the complex structure on $X$ and deform the K\"ahler class
to infinity. Let $f: X\rightarrow \mathbf{CP}^1$ be an elliptic
$K3$ surface. Let $\chi\geq0$ be the pullback of a K\"ahler form
on $\mathbf{CP}^1$ and $\omega_0$ be a K\"ahler form on $X$. We
construct a reference K\"ahler metric $\omega_t=\chi+t\omega_1$
and $[\omega_t]$ tends to $[\chi]$ as $t\rightarrow 0$. We can
always scale $\omega_1$ so that the volume of each fibre of $f$
with respect to $\omega_t$ is $t$. Suppose $\Omega$ is a
Ricci-flat volume form on $X$ with $\partial\dbar\log\Omega=0$.
Then Yau's proof \cite{Ya2} of Calabi's conjecture yields a unique
solution $\varphi_t$ to the following Monge-Amp\`ere equation for
$t\in (0, 1]$
\begin{equation}\label{K3}
\left\{
\begin{array}{rcl}
&& \frac{(\omega_t+\sqrt{-1}\partial\dbar\varphi_t)^2}{\Omega}=C_t\\
&&\int_{X}\varphi_t\Omega=0 , \end{array} \right.
\end{equation}
where $C_t=[\omega_t]^2$. Therefore we obtain a family of
Ricci-flat metrics
$\omega(t, \cdot)=\omega_t+\sqrt{-1}\partial\dbar\varphi_t$.
The following theorem is the main result of this section.
\begin{theorem}\label{structure}
Let $f: X\rightarrow \mathbf{CP}^1$ be an elliptically fibred $K3$
surface with $24$ singular fibres of type $I_1$. Then the
Ricci-flat metrics $\omega(t,\cdot)$ converges to the pullback of
a K\"ahler metric $\tilde{\omega}$ on $\mathbf{CP}^1$  in
any compact set of $X_{reg}$ in $C^{1,1}$ as $t\rightarrow 0$. The
K\"ahler metric $\tilde{\omega}$ on $\mathbf{CP}^1$ satisfies the
equation
\begin{equation}\label{K3eqn}
Ric(\tilde{\omega})=\omega_{WP}.
\end{equation}
\end{theorem}
\begin{proof} All the estimates can be obtained by the same argument in Section 4
with little modification. It is relatively easy compared to the
K\"ahler-Ricci flow because there is no such a term of
$\ddt{\varphi}$. We apply similar argument Section 6.1 to prove
the weak convergence. It is not difficult to show that  for any
test function $\zeta\in C^{\infty}_0(\mathbf{CP}^1_{reg})$  we
have
$$\int_{X}2\zeta (\chi+\sqrt{-1}\partial\dbar\varphi_t)\wedge\omega_{SF}
=\int_{X}\zeta\Omega +O(t).$$
By taking the convergent
subsequent we have
$$\frac{\chi+\sqrt{-1}\partial\dbar\varphi_0}{\chi}
=\frac{\Omega}{2\chi\wedge\omega_{SF}}.$$
Since such $\varphi_0$
with bounded $\|\varphi_0\|_{L^{\infty}}$ is unique and we have
the 2nd order estimates for $\varphi_t$, the convergence is then
uniform and this completes our proof. \qed
\end{proof}
This limit metric $\tilde{\omega}$ coincides with McLean's metric
as obtained by Gross and Wilson \cite{GrWi}.  Their construction
is certainly much more delicate and gives an accurate
approximation near the singular fibres by the Ooguri-Vafa metrics.
Also Mclean's metric is an example of the generalized
K\"ahler-Einstein metric defined in Definition \ref{KEmetric}
satisfying
$$Ric(\omega)=-\lambda \omega +\omega_{WP}$$
when $\lambda=0$.
%
%

\section{Generalizations and problems}

\subsection{A metric classification for surfaces of
non-negative Kodaira dimension}

In this section we will give a metric classification for surfaces
of non-negative Kodaira dimension.
Any surface $X$ with $K_X\geq 0$  must be minimal and $\nu(X)>0$.
\begin{enumerate}

\item When $\nu(X)=2$, $X$ is a minimal surface of general type
and we have the following theorem.

\begin{theorem}\label{tianzhou}\textnormal{\cite{TiZha}}
If $X$ is a minimal complex surface of general type, then the
global solution of the K\"ahler-Ricci flow converges to a positive
current $\omega_{\infty}$ which descends to the K\"ahler-Einstein
orbifold metric on its canonical model. In particular,
$\omega_{\infty}$ is smooth outside finitely many rational curves
and has local continuous potential.
\end{theorem}

\item When $\nu(X)=1$, $X$ is a minimal elliptic surface. Theorem
\ref{main} shows that the K\"ahler-Ricci flow deforms any K\"ahler
metric to a unique generalized K\"ahler-Einstein metric $\omega$
on its canonical model $X_{can}$.

\item When $\nu(X)=0$, $X$ is a Calabi-Yau surface. The
K\"ahler-Ricci defined in \cite{Ca} deforms any K\"ahler metric to
a Calabi-Yau metric in the same K\"ahler class.
\end{enumerate}
\ When $X$ is not minimal, the K\"ahler-Ricci flow (\ref{krflow1})
will develop singularities at finite time. Let $\omega_0$ be the
initial K\"ahler metric and $T$ be the first time such that
$e^{-t}[\omega_0]-(1-e^{-t})2\pi c_1(X)$ fails to be a K\"ahler
class. The K\"ahler-Ricci flow has a smooth solution $\omega(t,
\cdot)$ on $[0, T)$ converging to a degenerate metric as $t$ tends
to $T$ (cf. \cite{TiZha}, also see \cite{CaLa}). This degenerate
metric is actually smooth outside a subvariety $C$. Such a $C$ is
characterized by the condition that
$e^{-T}[\omega_0]-(1-e^{-T})2\pi c_1(X)$ vanishes along $C$. This
implies that $C$ is a disjoint union of finitely many rational
curves with self-intersection $-1$. Then we can blow down these
$(-1)$-curves and obtain a complex surface $X'$ and
$e^{-T}[\omega_0]-(1-e^{-T})2\pi c_1(X)$ descends to a K\"ahler
class on $X'$. Choose a K\"ahler metric $\omega_{T}$ representing
this class and then $(1-e^{-T})^{-1}(\omega_T-e^{-T}\omega_0)$
represents $-2\pi c_1(X)$. Define
$$\omega_t(\cdot)=\frac{e^{-t}-e^{-T}}{1-e^{-T}}\omega_0+\frac{1-e^{-t
}}{1-e^{-T}} \omega_T$$ and write $\omega(t,
\cdot)=\omega_t(\cdot)+\sqrt{-1}\partial\dbar \varphi(t,\cdot)$,
then $\varphi$ solves (2.8) on $[0, T)$, moreover, $\varphi$
converges to a bounded function $\varphi(T,\cdot)$ as $t$ tends to
$T$ and $\varphi(T)$ is smooth outside $C$. Applying the above
vanishing property of $C$ (cf. \cite{TiZha}), one can show that
$\varphi(T,\cdot)$ descends to a continuous function on $X'$. We
believe that $\varphi(T,\cdot)$ is actually $C^{1,1}$. Since
$\omega_T$ is a smooth metric on $X'$, we can consider the flow
(2.8) on $X'$ with $\varphi(T,\cdot)$ as the initial potential. We
expect that (2.8) still has a unique solution on $(0, T')$ where
$T'$ is the first time such that $e^{-t}[\omega_T]-(1-e^{-t})2\pi
c_1(X')$ fails to be a K\"ahler class, hence, (\ref{krflow1}) has
a smooth solution $\omega(t, \cdot)$ on $X'\times (0, T')$. Either
$T'$ is $\infty$ or we can repeat the previous procedure and
continue the flow (\ref{krflow1}). After finitely many times, we
will get a minimal complex surface with non-negative Kodaira
dimension. Then the flow has a global solution which falls into on
the cases described above. In order to complete this, one needs to
do the following steps: 1. prove an optimal estimate for
$\varphi(T, \cdot)$; 2. prove that the flow (2.8) has a unique
solution under weaker assumptions on smoothness of initial data
$\varphi(T, \cdot)$. Of course, if one can get sufficient
regularity in step 1, step 2 follows from the standard theory of
parabolic equations. We will address these problems in a
forthcoming paper.

\subsection{Higher dimension}

In this section, we discuss possible generalizations of Theorem
\ref{main} for higher dimension. First, as we assumed in Section
5.1, let $X$ be a non-singular variety of $\dim X=n$  such that
$\mu K_X$ is base point free for $\mu$ sufficiently large. Then
the pluricanonical map defines a holomorphic fibration $f:
X\rightarrow X_{can}$ by the linear system $|\mu K_X|$, where
$X_{can}$ is the canonical model of $X$.

\begin{enumerate}
\item If $\nu(X)=n$, $K_X$ is big and nef. Hence $X$ is a minimal
variety of general type. The K\"ahler-Ricci flow will deform any
K\"ahler metric to a canonical K\"ahler-Einstein metric on $X$
(cf. \cite{ TiZha, Ts}).

\item If $\nu(X)=1$, $X_{can}$ is a curve. With  little
modification of the proof, Theorem \ref{main} can be generalized
and the K\"ahler-Ricci flow will converge.

\item If $1<\mu(X)<n$, the fibration structure of $f$ can be very
complicated.  A large number of the calculations can be carried
out as in this paper and we expect the K\"ahler-Ricci flow will
converge appropriately to the pullback of a canonical metric
$\omega_{\infty}$ on the
 $X_{can}$ such that
$Ric(\omega_{\infty})=-\omega_{\infty}+\omega_{WP}$ on
$X^0_{can}$.
\end{enumerate}
When $K_X$ is nef, the K\"ahler-Ricci flow has long time
existence, yet it does not necessarily converge, although the
abundance conjecture predicts that $\mu K_X$ is globally generated
for $\mu$ sufficiently large. Hence, the problem of convergence of
the K\"hler-Ricci flow for nef $K_X$ can be considered as the
analytic version of the abundance conjecture.
If $K_X$ is not nef, the flow will develop finite time
singularities. Let $\omega_0$ be the initial K\"ahler metric and
$T$ the first time such that $e^{-t}[\omega_0]-(1-e^{-t})2\pi c_1(X)$
fails to be a K\"ahler class. The potential $\varphi(T,\cdot)$ is
bounded and smooth outside an analytic set of $X$(cf.
\cite{TiZha}). Let $X_1$ be the metric completion of $\omega(T,
\cdot)$. We conjecture that this is an analytic variety. It might
be a flip of $X$, or a variety obtained by certain standard
algebraic procedure. Of course $X_1$ might have singularities and
it is not clear at all how to develop the notion of a weak Ricci
flow on a singular variety. Suppose such a procedure can be
achieved and the K\"ahler-Ricci flow can continue on $X_1$, then
after applying the above procedure finitely many times on $X_1$,
$X_2$, ..., $X_N$, we might have $K_{X_N}\geq 0$ and get the
minimal model of $X$.

\bigskip
\noindent {\bf Acknowledgements} The first named author would like
to thank Professor D.H. Phong, his thesis advisor, for his
continued support, encouragement and advice, H. Fang, J. Sturm and
B. Weinkove for some enlightening discussions. He would also like
to thank the members of the complex and algebraic geometry group
at Johns Hopkins for their support.


\bigskip
\bigskip
\bigskip
\bigskip
\begin{thebibliography}{99}

\bibitem[Au]{Au} Aubin, T. {\em Equations du type Monge-Amp\`ere sur les vari\'et\'es K\"ahleriennes compacts}, Bull. Sc. Math. 102 (1976),
119-121.

\bibitem[BaMu]{BaMu} Bando, S. and Mabuchi, T., {\em Uniqueness of Einstein
K\"ahler metrics modulo connected group actions}, Algebraic
geometry, Sendai, 1985, 11--40, Adv. Stud. Pure Math., 10,
North-Holland, Amsterdam, 1987.

\bibitem[BaHuPeVa]{Ba} Barth, W., Hulek, K., Peters, C.,
and Van De Ven, A., {\em Compact complex surfaces}, 2003,
Springer.

\bibitem[Ca]{Ca}
Cao, H., {\em Deformation of K\"ahler metrics to K\"ahler-Einstein
metrics on compact K\"ahler manifolds}, Invent. Math. 81 (1985),
no. 2, 359--372.

\bibitem[CaLa]{CaLa} Cascini, P. and La Nave, P.,
{\em K\"ahler-Ricci Flow and the Minimal Model Program for
Projective Varieties}, preprint.

\bibitem[Ch]{Ch} Chow, B., {\em The Ricci flow on the $2$-sphere},
 J. Differential
Geom. 33 (1991), no. 2, 325--334.

\bibitem[ChTi]{ChTi} Chen, X.X. and Tian, G.,
{\em Ricci flow on K\"ahler-Einstein surfaces}, Invent. Math. 147
(2002), no. 3, 487--544.

\bibitem[ChYa]{ChYa} Cheng, S. Y. and Yau, S. T.,
{\em Differential equations on Riemannian
manifolds and their geometric applications}, Comm. Pure Appl.
Math. 28 (1975), no. 3, 333--354.

\bibitem[ClKoMo]{ClKoMo} Clemens, H,, Kollar, J. and Mori, S., {\em
Higher-dimensional complex geometry}, Ast\'eisque No. 166 (1988),
144 pp. (1989).

\bibitem[Do]{Do} Donaldson, S. K., {\em Scalar curvature and projective
embeddings, I.}, J. Differential Geom. 59 (2001), no. 3, 479-522.

\bibitem[GrWi]{GrWi} Gross, M. and Wilson, P. M. H., {\em Large complex structure limits of
$K3$ surfaces}, J. Differential Geom. 55 (2000), no. 3, 475--546.

\bibitem[Ha1]{Ha1}
Hamilton, R., {\em Three-manifolds with positive Ricci curvature},
J. Differential Geom. 17 (1982), no. 2, 255--306.

\bibitem[Ha2]{Ha2}
Hamilton, R., {\em  The formation of singularities in the Ricci
flow}, Surveys in differential geometry, Vol. II (Cambridge, MA,
1993), 7--136, Internat. Press, Cambridge, MA, 1995.

\bibitem[H\"o]{Ho} H\"ormander, L., {\em An introduction to complex analysis
in several variables}, Van Nostrand, Princeton 1973.

\bibitem[Ko]{Ko} Kodaira, K., { On compact complex analytic surfaces I}, Ann. of
Math. (2) 71 (1960), 111--152. {\em II},  Ann. of Math. (2) 77
(1963), 563--626. {\em III}, Ann. of Math. (2) 78 (1963), 1--40.

\bibitem[Kol1]{Ko1} Kolodziej, S., {\em The complex Monge-Amp\`ere
equation}, Acta Math. 180 (1998), no. 1, 69--117.

\bibitem[Kol2]{Ko2} Kolodziej, S., {\em Stability of solutions to
the complex Monge-Amp\`ere equation on compact K\"ahler manfiolds},
preprint.

\bibitem[KoSo]{KoSo} Kontsevich, M. and Soibelman, Y., {\em Homological mirror symmetry and torus fibrations}, Symplectic geometry and
mirror symmetry (Seoul, 2000), 203--263, World Sci. Publishing,
River Edge, NJ, 2001.

\bibitem[LiYa]{LiYa} Li, P. and Yau, S.T., {\em Estimates of eigenvalues of a compact Riemannian manifold}, Geometry of the Laplace operator
(Proc. Sympos. Pure Math., Univ. Hawaii, Honolulu, Hawaii, 1979),
pp. 205--239, Proc. Sympos. Pure Math., XXXVI, Amer. Math. Soc.,
Providence, R.I., 1980.

\bibitem[Mi]{Mi} Miranda, R., {\em The basic theory of elliptic
surfaces}, Dottorato di Ricerca in Matematica. [Doctorate in
Mathematical Research] ETS Editrice, Pisa, 1989. vi+108 pp.

\bibitem[Mc]{Mc} McLean, R., {\em Deformations of calibrated
submanifolds}, Comm. Anal. Geom. 6 (1998), no. 4, 705--747.

\bibitem[OoVa]{OoVa} Ooguri, H. and Vafa, C., {\em Summing up Dirichlet instantons}, Phys. Rev. Lett. 77 (1996), no. 16, 3296--3298.

\bibitem[Pe]{Pe} Perelman, P., {\em The entropy formula for the Ricci flow
and its geometric
applications}, preprint math.DG/0211159.
\bibitem[PhSt]{PhSt} Phong, D. H. and Sturm, J., {\em On stability and the
convergence of the K\"ahler-Ricci flow
}, preprint math.DG/0412185.

\bibitem[SeTi]{SeTi} Sesum, N. and Tian, G., {\em Bounding scalar
curvature and diameter along the K\"ahler Ricci flow (after
Perelman)}, lecture note.

\bibitem[Si]{Si} Siu, Y-T. {\em Lectures on Hermitian-Einstein metrics for
stable  bundles and K\"ahler-Einstein metrics}, Birkh\"auser
Verlag, Basel 1987.

\bibitem[SoWe]{SoWe} Song, J. and Weinkove, B., {\em On
the convergence and singularities of the J-flow with applications to the Mabuchi
energy}, preprint math.DG/0410418.

\bibitem[StYaZa]{SYZ} Strominger, A., Yau, S.T. and Zaslow, E., {\em Mirror
symmetry is $T$-duality}, Nuclear Phys. B 479 (1996), no. 1-2,
243--259.
\bibitem[Ti1]{Ti1} Tian, G., {\em On K\"ahler-Einstein metrics on certain K\"ahler
manifolds with $C_1(M)>0$}, Invent. Math. 89 (1987), no. 2,
225--246.

\bibitem[Ti2]{Ti2} Tian, G., {\em On Calabi's conjecture for complex surfaces
with positive first Chern class}, Invent. Math. 101, no. 1 (1990),
101-172.

\bibitem[TiZha]{TiZha} Tian, G. and Zhang, Z., {\em A note on
the K\"ahler-Ricci flow on
projective manifolds of general type}, preprint.

\bibitem[TiZhu]{TiZhu} Tian, G. and Zhu, X., {\em Convergence of
K\"ahler Ricci flow}, preprint.

\bibitem[Ts]{Ts} Tsuji, H. {\em Existence and degeneration of K\"ahler-Einstein metrics on minimal algebraic varieties of general type}, Math.
Ann. 281 (1988), 123-133.

\bibitem[Ya1]{Ya1}Yau, S.T., {\em A general Schwarz lemma for
K\"ahler manifolds}, Amer. J. Math. 100 (1978), no. 1, 197--203.

\bibitem[Ya2]{Ya2} Yau, S.T. {\em On the Ricci curvature of a compact K\"ahler manifold and the complex Monge-Amp\`ere equation, I}, Comm. Pure Appl. Math. 31 (1978), 339-411.

\bibitem[Ya3]{Ya3} Yau, S.T., {\em Open problems in geometry}, Proc. Symposia Pure
Math. 54 (1993), 1-28.

\bibitem[Zh]{Zh} Zhang, Z., Thesis at MIT.

\end{thebibliography}
\end{document}